\documentclass[a4paper,12pt]{article}
\usepackage{amssymb}
\usepackage{amsmath} 
\usepackage{amsthm}
\usepackage{makeidx}
\usepackage{amscd}

\font\rus=wncyr10 scaled \magstep1%

\def\hfl#1#2#3{\smash{\mathop{\hbox
to#3{\rightarrowfill}}\limits
^{\scriptstyle#1}_{\scriptstyle#2}}}

\def\vfl#1#2#3{\llap{$\scriptstyle #1$}
\left\downarrow\vbox to#3{}\right.\rlap{$\scriptstyle #2$}}

\theoremstyle{plain}

\newtheorem{theo}{Theorem}[section]

\newtheorem{prop}[theo]{Proposition}
\newtheorem{lem}[theo]{Lemma}
\newtheorem{cor}[theo]{Corollary}

\newtheorem{defi}[theo]{Definition}

\renewcommand{\thesection}{\arabic{section}.}

\newtheorem{theorem}{Theorem}[section]
\newtheorem{proposition}[theorem]{Proposition}
\newtheorem{lemma}[theorem]{Lemma}

\theoremstyle{definition}

\usepackage{verbatim}
\def \kbar {{\overline{k}}}

\def \Romannumeral #1 {\expandafter\uppercase\expandafter {\romannumeral #1} }

\def \Br {{\rm{Br\,}}}

\def \Ga {{\Gamma}}  
\def \R {{\bf R}}
\def \Pic {{\rm {Pic\,}}}
\def \Div {{\rm {Div\,}}}

\def \Gal {{\rm{Gal\,}}}
\def \Ker {{\rm{Ker\,}}}

\def \A{{\mathbb A}}
\def \St {{\rm{St}}}
\def \tor {{\rm{tor}}}

\def \Hom {{\rm {Hom}}}

\def \SL {{\rm {SL}}}

\def \PGL {{\rm {PGL}}}

\def\ov{\overline}

\def\k{\ov k}
\def \Z {{\bf Z}}
\def \Q {{\bf Q}}

\def \Br {{\rm Br\,}}
\def \Ext {{\rm Ext}}

\def\G{{\bf G}}
\def\R{{\bf R}}
\def\P{{\bf P}}

\def\nr{{\textup{nr}}}

\def\H{{\rm H}}

\def\Sha{{\cyr X}}

\def\ab{{\rm ab}}
\def\der{{\rm der}}
\def\mult{{\rm mult}}

\def\Ga{{\mathfrak{g}}}

\newcommand{\bthe}{\begin{theo}}
\newcommand{\ble}{\begin{lem}}
\newcommand{\bpr}{\begin{prop}}
\newcommand{\bco}{\begin{cor}}
\newcommand{\bde}{\begin{defi}}
\newcommand{\ethe}{\end{theo}}
\newcommand{\ele}{\end{lem}}
\newcommand{\epr}{\end{prop}}
\newcommand{\eco}{\end{cor}}
\newcommand{\ede}{\end{defi}}

\newcommand{\mathb}[1]{\mathbf{#1}}

\newcommand{\capbar}{\overline}

\newcommand{\sV}{\mathcal{V}}
\newcommand{\sY}{\mathcal{Y}}
\newcommand{\isoto}{\overset{\sim}{\to}}

\newcommand{\into}{\hookrightarrow}

\newcommand{\Zz}{{\mathb{Z}}}
\newcommand{\Qq}{{\mathb{Q}}}
\newcommand{\Rr}{{\mathb{R}}}
\newcommand{\Cc}{{\mathb{C}}}

\newcommand{\Gm}{\mathb{G}_m}

\def\uu{^\mathrm{u}}

\def\red{^\mathrm{red}}
\def\tor{^{\mathrm{tor}}}

\def\sss{^{\mathrm{ss}}}

\def\ssu{^{\mathrm{ssu}}}

\def\sa{^{\mathrm{sab}}}
\def\sab{^{\mathrm{sab}}}
\def\lin{^{\mathrm{lin}}}
\DeclareMathOperator{\NS}{NS}

\newcommand{\et}{\etale}

\newcommand{\Zbar}{{\capbar{Z}}}

\newcommand{\Gbar}{{\capbar{G}}}
\newcommand{\Hbar}{{\capbar{H}}}

\def\et{{\textup{\'et}}}
\def\Zar{{\textup{Zar}}}

\def\xbar{{\overline{x}}}

\def\pibar{\overline{\pi}}
\def\abar{\overline{a}}

\renewcommand{\infty}{{r}}

\def\NS{\textup{NS\,}}
\def\Bro{\textup{Br}_1}

\def\bsm{\left( \begin{smallmatrix}}
\def\esm{\end{smallmatrix} \right)}

\makeindex

\DeclareFontFamily{U}{wncy}{}
\DeclareFontShape{U}{wncy}{m}{n}{%
<5>wncyr5%
<6>wncyr6%
<7>wncyr7%
<8>wncyr8%
<9>wncyr9%
<10>wncyr10%
<11>wncyr10%
<12>wncyr6%
<14>wncyr7%
<17>wncyr8%
<20>wncyr10%
<25>wncyr10}{}
\DeclareMathAlphabet{\cyr}{U}{wncy}{m}{n}

\title{The elementary obstruction and homogeneous spaces}
\author{M. Borovoi, J-L. Colliot-Th\'el\`ene and A.N. Skorobogatov}
\date{}
\begin{document}
\baselineskip=15pt
\maketitle

MSC Primary: 14G05, 11G99, 12G05

MSC Secondary: 11E72, 14F22, 14K15, 14M17, 20G99

\medskip

{\bf Abstract}

Let $k$ be a field of characteristic zero and
  ${\overline  k}$   an algebraic closure of  $k$.
For a  geometrically integral variety $X$ over $k$,
we write ${\overline  k}(X)$ for the function field of ${\overline  X}=X\times_k{\overline  k}$.
If $X$ has a smooth $k$-point, the
 natural embedding of multiplicative groups
${\overline  k}^*\hookrightarrow  {\overline  k}(X)^* $ admits a Galois-equivariant retraction.

In the first part of the paper, over  local  and then over global fields,
equivalent conditions to the existence of such  a retraction are given.
They are expressed in terms of the Brauer group of $X$.

In the second part of the paper, we restrict attention
to varieties   which are homogeneous spaces of connected but otherwise arbitrary algebraic groups,
 with connected
geometric stabilizers. For $k$ local or global,
for such a variety $X$,  in many situations but not all, the existence of
a Galois-equivariant retraction to  ${\overline  k}^*\hookrightarrow {\overline  k}(X)^* $ ensures the existence of a $k$-rational point on $X$.
For homogeneous spaces of linear algebraic groups,
the technique also handles the case where  $k$ is the function field of a complex surface.

\bigskip
{\bf R\'esum\'e}

Soient $k$ un corps de caract\'eristique nulle et ${\overline  k}$   une cl\^oture alg\'ebrique de  $k$.
Pour une $k$-vari\'et\'e  $X$ g\'eom\'etriquement int\`egre,  on note
${\overline  k}(X)$ le corps des fonctions de ${\overline  X}=X\times_k{\overline  k}$.
Si $X$ poss\`ede un $k$-point lisse, le plongement naturel de groupes multiplicatifs
${\overline  k}^*\hookrightarrow {\overline  k}(X)^* $ admet une r\'etraction \'equivariante pour l'action du groupe de Galois de ${\overline  k}$
sur $k$.

Dans la premi\`ere partie de l'article, sur les corps locaux puis sur les corps globaux,
on donne des conditions \'equivalentes \`a l'existence d'une telle r\'etraction \'equivariante.
Ces conditions s'expriment en terme du groupe de Brauer de la vari\'et\'e $X$.

Dans la seconde partie de l'article, on consid\`ere le cas des espaces homog\`enes de groupes alg\'ebriques connexes, non n\'ecessairement lin\'eaires,
avec groupes d'isotropie g\'eom\'etriques connexes. Pour $k$ local ou global, pour un tel espace homog\`ene $X$, dans beaucoup de cas
mais pas dans tous, l'existence d'une r\'etraction \'equivariante \`a ${\overline  k}^*\hookrightarrow  {\overline  k}(X)^* $  implique l'existence
d'un point $k$-rationnel sur $X$.  Pour les espaces homog\`enes de groupes lin\'eaires, la technique permet aussi de traiter le cas o\`u $k$
est un corps de fonctions de deux variables sur les complexes.

\section{Introduction}

Among the many obstructions to the existence of rational points
one is particularly remarkable by the simplicity of construction.

Let $k$ be a field of characteristic zero,
  $\ov k$   an algebraic closure of  $k$ and
$\Ga$  the Galois group of $\ov k$ over $k$.
For a  geometrically integral variety $X$ over $k$,
we write $\ov k (X)$ for the function field of $\ov X=X\times_k\ov k$.
The elementary obstruction, defined and studied in \cite{CS}, is
the class $ob(X)\in\Ext^1_\Ga(\kbar(X)^*/\kbar^*,\kbar^*)$
of the extension of Galois modules
\begin{equation}
1\to\ov k^*\to \ov k(X)^*\to \ov k(X)^*/\ov k^*\to 1. \label{ob}
\end{equation}

 If $X$ and $Y$ are geometrically integral $k$-varieties
and there exists a dominant rational map $f$ from $X$ to $Y$, then $ob(X)=0$
implies  $ob(Y)=0$. In particular, the vanishing of $ob(X)$ is
a birational invariant of $X$.
As pointed out by O. Wittenberg (\cite{Witt},  Lemma 3.1.2), one can say more:  if there exists a rational map from a
geometrically integral variety $X$ to a smooth geometrically integral
$k$-variety $Y$, then $ob(X)=0$ implies  $ob(Y)=0$.

 As a special case, if $X$ has a smooth $k$-point, the extension (\ref{ob})  is split, so that
$ob(X)=0$ (\cite{CS}, Prop. 2.2.2).

We are thus confronted with
the following natural question:
for which fields $k$ and $k$-varieties $X$
is $ob(X)$ the only obstruction to the existence of $k$-points on $X$?

In the first part of this paper we consider arbitrary smooth,
geometrically integral varieties.
After recalling some general facts about the elementary obstruction,
we turn to local and global fields. For such fields we relate
 the elementary obstruction to
obstructions coming from the Brauer group:

(i) If $k$ is local (e.g., a $p$-adic field or the   field of real numbers),
$ob(X)=0$ if and only if the natural map $\Br k \to \Br k(X)$ is injective
(see Theorems \ref{t} and \ref{treal} for more general statements).

(ii) If $k$ is a number field,  if $ob(X)=0$
and $X$ has points in all completions of $k$, then any ad\`ele of $X$
is orthogonal
to the subgroup of the Brauer group of $X$ consisting of  `algebraic' elements
which are everywhere locally constant (Theorem \ref{B}).

In the second part of the paper we explore the elementary
obstruction $ob(X)$, when $X$ is a
homogeneous space of a connected algebraic
group $G$, not necessarily linear. Most results require the assumption
that the stabilizers of $\ov k$-points of $X$ are connected.
Under this assumption, we prove the following results.

(iii) If $k$ is a $p$-adic field, we show that $ob(X)=0$ implies the existence of
a rational point (Theorem \ref{t2}). This actually holds as long as the Brauer group of $k$ injects into
the Brauer group of the function field of $X$ (Corollary \ref{cor2}). The case of
homogeneous spaces of abelian varieties was known
(Lichtenbaum,   van Hamel).

(iv) If $k$ a `good'  field of cohomological dimension at most 2,
and the group $G$ is  linear,
the hypothesis $ob(X)=0$ implies the existence of a rational point
(Theorem \ref{good}).
This result covers the case of $p$-adic fields (already handled in (iii))
and of totally imaginary number fields.
Thanks to a theorem of de Jong, it also applies to function fields in two variables
over an algebraically closed field, provided that $G$ has no factor of type $E_{8}$.

(v) If  $k$ is a number field and the group $G$ is linear,
 if $X$ has points in the real completions of $k$ and $ob(X)=0$, then
$X$ has a rational point (Theorem \ref{real}).

(vi) If $k$ is a totally imaginary number field and
$G$ is an arbitrary connected algebraic group,
assuming finiteness of the Tate--Shafarevich group of
the maximal abelian variety quotient of $G$, we prove
that $ob(X)=0$ implies that $X$
has a rational point (Theorem \ref{imago}).
A key ingredient is a recent result of
D. Harari and T. Szamuely on principal homogeneous spaces of
commutative algebraic groups. Their theorem also holds when
$k$ has real completions.

(vii) In the general case of arbitrary connected groups we
found, somewhat to our surprise, a principal homogeneous
space $X/\Q$ of a non-commutative group $G$
with $ob(X)=0$, with points everywhere locally, but without
$\Q$-points (Proposition \ref{example}). By a previously mentioned result
(Theorem \ref{B}, or by an easy direct argument) one sees that
the Brauer--Manin obstruction
attached to the subgroup ${\cyr B}(X)\subset \Bro X$
of everywhere locally constant classes, is trivial. Thus we obtain a negative
answer to the following question raised in \cite{S} (Question 1, p. 133):
is the Brauer--Manin obstruction attached to ${\cyr B}(X)$ the only
obstruction to the Hasse principle for torsors of
arbitrary connected algebraic groups?
This phenomenon is due to a combination of three factors:
the presence of real places, the
non-commutativity and the non-linearity of $G$.

The example in (vii) can be accounted for by the Brauer--Manin obstruction
attached to the group $\Bro X^c$, where $X^c$ denotes a smooth compactification
of the torsor $X$. This is a special instance of a result of Harari \cite{H}.
In the Appendix,
building upon \cite{H} and the techniques
in the present paper, we extend Harari's result to  homogeneous spaces
of  any connected algebraic group $G$, assuming that the geometric stabilizers
are connected. As in  \cite{H} and earlier work on the subject, the result
here is conditional on the  finiteness of the Tate--Shafarevich group of
the maximal abelian variety quotient of $G$.

\medskip

In the case of a linear algebraic group $G$, the recurring assumption that
the geometric stabilizers are connected can be somewhat relaxed
(Theorems \ref{thm:non-connected} and \ref{thm:non-connected-Manin}), but some condition must definitely be imposed,
as shown by an example of M. Florence \cite{Fl}.

The starting point of our work was the following
result of Joost van Hamel:
for a principal homogeneous space $X$ of a connected linear $k$-group $G$
over a $p$-adic field $k$, the elementary obstruction is the only obstruction
to the existence of a $k$-rational point in $X$.

The authors are grateful to Mathieu Florence, Joost van Hamel,
David Harari, Olivier Wittenberg and Yuri Zarhin for useful discussions.
We thank the Mathematical Sciences Research Institute at Berkeley,
where this work was done, for hospitality
and support.
The first named author was partially supported by the Hermann Minkowski Center for Geometry.
The second named author acknowledges the support of the Clay Institute.

\section{Elementary obstruction}

\subsection{Preliminaries}

Let $k$ be a field of characteristic 0, $\ov k$
an algebraic closure of $k$, $\Ga=\Gal(\ov k/k)$.
If $X$ is a $k$-variety, we let $\ov X = X \times_{k}\ov k$.
If $X$ is integral, we denote by $k(X)$ the function field of $X$.
If $X$ is geometrically integral, we denote by $\ov k (X)$ the function field
of $\ov X$.
We let $\Div X$ denote the group of Cartier divisors on $X$,
$\Pic X$ denote the Picard group $\H^1_\Zar(X,\G_{m})= \H^1_{\et}(X,\G_{m})$  of $X$. By $\Br X$ we denote
the cohomological Brauer--Grothendieck group $\H^2_{\et}(X,\G_{m})$,
and by $\Bro X$ the kernel of the natural map
$\Br X \to \Br \ov X$. If $M$ is a continuous discrete $\Ga$-module,
we write $\H^i(k,M)$ for the Galois cohomology groups.

When $\ov k^* = \ov k[X]^*$ the Hochschild--Serre
spectral sequence
$$E_{2}^{pq}=\H^p(k,\H_{\et}^q(\ov X,\G_m))\Rightarrow
\H_{\et}^{p+q}(X,\G_m)$$
gives rise to the well known exact sequence
\begin{equation}
0\to\Pic X\to (\Pic\ov X)^\Ga\to\Br k\to \Bro X
\hfl{r}{}{6mm}\H^1(k,\Pic\ov X), \label{e1}
\end{equation}
where the  map $\Bro  X \to \H^1(k,\Pic\ov X)$
is onto if $X$ has a $k$-point, or
if  $k$ is a local or global field.

Recall that if $A$ and $B$ are continuous discrete $\Ga$-modules,
then $\Ext^n_\Ga(A,B)$ is defined as the derived functor
of $\Hom_\Ga(A,B)$
in the second variable (see  \cite{EC}).
In particular, there are long exact sequences in either variable,
and the elements of $\Ext^n_\Ga(A,B)$ classify equivalence
classes of $n$-extensions of continuous discrete Galois modules.

\medskip

Let $X$ be a smooth, quasi-projective and
geometrically integral variety over $k$.
Then Cartier divisors coincide with Weil divisors, which implies that
$\Div\ov X$ is a permutation $\Ga$-module.
We have the following natural 2-extension of
continuous discrete $\Ga$-modules:
$$
1\to \ov k[X]^*\to \ov k(X)^*\to \Div\ov X\to\Pic\ov X\to 0.
$$
When $\ov k^* = \ov k[X]^*$, this reads
\begin{equation}
1\to\ov k^*\to \ov k(X)^*\to \Div\ov X\to\Pic\ov X\to 0.\label{e(X)}
\end{equation}
Under the assumption $\ov k^* = \ov k[X]^*$,
write $e(X)\in \Ext^2_\Ga(\Pic\ov X,\ov k^*)$ for
the corresponding class\footnote{This
definition of $e(X)$ differs from that in \cite{S} by $-1$.}.
Much is known about the classes $ob(X)$ and $e(X)$ (see \cite{CS}, Section 2,
or \cite{S}, Ch. 2). Clearly $e(X)$ is the cup-product of
$$
1\to \ov k(X)^*/k^* \to \Div\ov X\to\Pic\ov X\to 0$$
with  the class $ob(X)$.
For further reference we list here
some of the known properties of these classes.

\ble \label{facts1}

{\rm (i)}
The class $ob(X)$ lies in the kernel of the
natural
map
$$\Ext^1_\Ga(\ov k(X)^*/\ov k^*,\ov k^*) \to \Ext^1_\Ga(\ov k(X)^*/\ov k^*,{\ov k}(X)^*).$$

{\rm (ii)} If there exists a zero-cycle of
degree $1$ on $X$, then $ob(X)=0$.

{\rm (iii)} If $ob(X)=0$, then
for a $k$-group of multiplicative type $S$ and $i=0,1,2$ the natural
maps $\H^i(k,S) \to \H^i(k(X),S)$ are injective. In particular,
 the map $\Br k \to \Br k(X)$ is injective, and so is the map
$\Br k \to \Br X$.

{\rm (iv)} If $X$ is $k$-birational to a homogeneous space of a $k$-torus,
then $ob(X)=0$ if and only if $X(k) \neq \emptyset$.

\ele

{\it Proof\/} (i) is obvious.

(ii)
\cite{CS}, Prop. 2.2.2     (see also \cite{S}, Thm. 2.3.4).

(iii) \cite{CS}, Prop. 2.2.5.

(iv) We may assume that $X$ is a $k$-torsor of a $k$-torus (cf. \cite{B96}, proof of Prop. 3.3).
If $ob(X)=0$, then $\ov k^*$ is a direct summand in $\ov k(X)^*$,
hence it is also a direct summand in $\ov k[X]^*$. Now it follows from
\cite{San}, (6.7.3) and (6.7.4), that
$X$ is a trivial torsor, i.e. $X$ has a $k$-point. QED

 \ble \label{facts2}
Assume $\ov k^* = \ov k[X]^*$.

{\rm (i)} $ob(X)=0$ if and only if $e(X)=0$.

{\rm (ii)} The class $e(X)$ lies in the kernel of the
natural
map
$$\Ext^2_\Ga(\Pic\ov X,\ov k^*) \to \Ext^2_\Ga(\Pic\ov X,{\ov k}(X)^*).$$

{\rm (iii)} The map $(\Pic\ov X)^\Ga\to\Br k$ in $(\ref{e1})$ is
the Yoneda product with $e(X)$ (up to sign).

{\rm (iv)} If $\Pic \ov X$ is finitely generated  and free as an abelian group, and $S$ denotes the $k$-torus with character group
$\Pic \ov X$, then $ob(X)=0$ if and only if $\H^2(k,S)$ injects into $\H^2(k(X),S)$.

{\rm (v)} If $\Pic\ov X=\Z$, then $ob(X)=0$ if and only if the map
$\Br k \to \Br k(X)$ is injective.

{\rm (vi)} If $\Pic \ov X$ is finitely generated and is a direct factor of a permutation
$\Ga$-module, then $ob(X)=0$ if and only if for any finite field extension $K/k$ the map
$\Br K \to \Br K(X)$ is injective.

{\rm (vii)} If $\Pic \ov X =0$, then
$ob(X)=0$.

{\rm (viii)} If $X$ is a principal homogeneous space of a
semisimple simply connected group, then $ob(X)=0$.

{\rm (ix)} If $X \subset \P^n_{k}$ is a smooth projective
hypersurface and $n \geq 4$, then $ob(X)=0$.

\ele

{\it Proof\/} (i)  \cite{CS}, Prop. 2.2.4, or \cite{S}, Thm. 2.3.4.

(ii) follows from (i) in the previous lemma.

(iii) \cite{CS}, Lemme 1.A.4, or \cite{S1}, Prop. 1.1.

(iv) The direct implication follows from (iii) of the previous lemma.
For the converse observe that
the natural map $\H^2(k,S)\to \H^2(k(X),S)$ factors through
\begin{equation}
\H^2(k,\Hom_\Z(\Pic\ov X,\ov k^*))\to\H^2(k,\Hom_\Z(\Pic\ov X,\ov k(X)^*)).
\label{hhh}
\end{equation}
Since the $\Ga$-module $\Pic\ov X$ is finitely generated
we have the spectral sequence
$$\H^p(k,\Ext^q_\Z(\Pic\ov X,\ov k(X)^*))\Rightarrow \Ext^{p+q}_\Ga(\Pic\ov X,\ov k(X)^*).$$
Since $\Pic\ov X$ is    finitely generated and torsion-free
we have $\Ext^q_\Z(\Pic\ov X,\ov k(X)^*)=0$ for any $q\geq 1$,
so that the spectral sequence degenerates and gives an isomorphism
$\H^2(k,\Hom_\Z(\Pic\ov X,\ov k(X)^*))=\Ext^2_\Ga(\Pic\ov X,{\ov k}(X)^*)$.
This, and a similar argument for $\H^2(k,\Hom_\Z(\Pic\ov X,\ov k^*))$,
identify (\ref{hhh}) with the map in (ii).
Now our statement follows from (i) and (ii).

(v) This is a special case of (iv).

(vi) Assume $ob(X)=0$. Let $K/k$ be a finite field extension. Applying
Lemma \ref{facts1} (iii) to the $k$-torus $S=R_{K/k}\G_{m}$ and using Shapiro's lemma,
one finds that  $\Br K \to \Br K(X)$ is injective. One could also directly argue that
$ob(X)=0$ implies $ob(X\times_{k}K)=0$.

 Assume now that  $\Pic \ov X$ is finitely generated and is a direct factor of a permutation
$\Ga$-module $\oplus_{i}\Z[\Ga/\Ga_{i}]$, where $\Ga_{i}={\rm Gal}(\ov k/K_{i})$, with each $K_{i} \subset \ov k$ a finite field extension of $k$. Let $S$, resp. $P$, be the $k$-torus whose character group is $\Pic \ov X$, resp. $\oplus_{i}\Z[\Ga/\Ga_{i}]$. There exists a $k$-torus $S_{1}$ and
an isomorphism of $k$-tori $S \times_{k} S_{1} \simeq P$.  Let us assume that for each $K_{i}/k$
the natural map $\Br K_{i} \to \Br K_{i}(X)$ is injective. By Shapiro's lemma this is equivalent
to assuming the injectivity of  the natural map $\H^2(k,P) \to \H^2(k(X),P)$. This in turn  implies the injectivity of
$\H^2(k,S) \to \H^2(k(X),S)$. From (iv) we conclude $ob(X)=0$.

(vii) Given (\ref{e(X)}),  this is an application of (i) (cf.
  \cite{CS},  Remarque 2.2.7.)

(viii) This is a direct application of (vii).

(ix) For such a hypersurface, the restriction map $\Z =\Pic\P^n_{k} \to \Pic X$
is an isomorphism, and so it is over $\ov k$ (Max Noether's theorem). Let $U \subset X$
be the complement of a smooth hyperplane section defined over $k$. Then $\ov k^*=\ov k [U]^*$ and $\Pic\ov U=0$.
One may then apply (vii) to $U$.
QED

\medskip

{\it Remarks}
(1) There exist higher Galois cohomological obstructions
to the existence of rational points,
and, more generally, to the existence of a zero-cycle of degree 1.
Let $X$ be a smooth geometrically integral $k$-variety,
and $S$ a $k$-group of multiplicative type, for instance
a finite $\Ga$-module.
 If $X$ has a zero-cycle of degree~1, then
for any positive integer $n$ the restriction map
$\H^n(k,S) \to \H^n(k(X),S)$ is injective:
this is a consequence of the Bloch--Ogus theorem.

(2)
In \cite{CS}, Exemples 2.2.12, one will find a sample of varieties,
over suitable fields, which satisfy $ob(X)=0$ but which have
no $k$-rational points.
Simple examples with $ob(X)=0$ are given by (viii) and (ix)
in the previous lemma. Some of these examples can
be explained by means of the higher Galois obstructions in Remark
(1), whereas some others cannot.
For more on this, see the Remarks after Theorems \ref{t} and \ref{treal}
below.

(3) Let $k={\bf C}((t))$. Let $X/k$ be the curve of genus 1 defined by
the homogeneous equation $x^3+ty^3+t^2z^3=0$. We obviously have $X(k)=\emptyset$.
The Brauer group of $k$ and of any finite extension
of $k$ vanishes. A general result of O. Wittenberg \cite{Witt} then ensures $ob(X)=0$.
Thus the absence of $k$-points on $X$ is not
detected by any of the above Galois cohomology arguments.

\medskip

{\it Questions}
 Let $X$ be a geometrically integral $k$-variety. Let $K/k$ be an arbitrary
field extension.

(1) Assume $ob(X)=0$. Does the $K$-variety $X_K$ satisfy $ob(X_K)=0$?
This is clear if $K\subset \ov k$.

(2) Assume that the $K$-variety $X_K$ satisfies $ob(X_K)=0$.
If the extension $K/k$ has
a $k$-place, does the $k$-variety $X$ satisfy $ob(X)=0$?

\medskip
We can answer the first question in a special case.
\bpr \label{obXK}
Let $X/k$ be a smooth, projective, geometrically integral variety. Assume that the Picard variety of $X$ is trivial.
If $ob(X)=0$ then for any field $K$ containing $k$ we have $ob(X_{K})=0$.
\epr

{\it Proof\/} Let $\k \subset \overline{K}$ be an inclusion of algebraic closures. Let ${\mathfrak{G}}={\rm Gal}(\overline{K}/K)$ and $ \Ga={\rm Gal}(\k/k)$. There is a natural map  ${\mathfrak{G}} \to \Ga$.
Because the Picard variety of $X$ is trivial, the  abelian groups $\Pic X_{\k} $ and  $\Pic X_{\overline{K}}$ are abelian groups of finite type and the
natural map
$\Pic X_{\k} \to \Pic X_{\overline{K}}$   is a Galois equivariant isomorphism  (the N\'eron--Severi group does not change under extensions of algebraically closed ground fields). There is an    equivariant commutative diagram of 2-extensions
\begin{equation}
\begin{array}{ccccccccccc}
1 &   \to &    \k^*              &   \to    &  \k(X)^*                      &         \to    &  \Div X_{\k}                    &        \to     &  \Pic  X_{\k}                   &   \to  &       0\\
   &          &\downarrow  &           &\downarrow               &                 &   \downarrow               &                      & \downarrow {\simeq}            &            &       \\
1&\to      & \overline{K}^*  &\to   &   \overline{K}(X)^*   &       \to      &  \Div X_{\overline{K}}   &   \to          &  \Pic X_{\overline{K}}    &   \to  & 0.
\end{array}
\label{obX_kK}
\end{equation}
If $ob(X)=0$ then the top 2-extension is trivial (Lemma \ref{facts2} (i)).
This implies that the bottom 2-extension is trivial, that is $ob(X_{K})=0$.
QED

Other  cases where Question (1) can be answered positively will be handled in  Subsections
\ref{sec:local-general} and \ref{sec:global-general} For further results, see \cite{Witt}.
\footnote{O. Wittenberg has very recently shown that the answer to Question (1), in general, is in the 
negative.}

\bigskip

Let $X/k$ be a smooth, projective, geometrically integral variety.
Let $J/k$ be the Picard variety of $X$.
Let $\NS \ov X$ be the
N\'eron--Severi group of $\ov X$.
{From} the exact sequence of $\Ga$-modules
\begin{equation}
0 \to J(\ov k) \to {\rm Pic\,}\ov X \to \NS \ov X \to 0 \label{1}
\end{equation}
we deduce the following diagram in which
the vertical sequences are exact:
\begin{equation}
\begin{array}{ccccc}
\H^1(k ,J)&\times & \Ext^1_\Ga(J(\ov k),\ov k^*)&\to &\Br k\\
\uparrow&&\downarrow&&||\\
(\NS \ov X)^\Ga &\times &\Ext^2_\Ga(\NS \ov X,\ov k^*)&\to &\Br k\\
\uparrow&&\downarrow&&||\\
({\rm Pic\,} \ov X)^\Ga &\times &\Ext^2_\Ga({\rm Pic\,} \ov X,\ov k^*)&\to &\Br k\\
\uparrow&&\downarrow&&||\\
J(k) &\times &\Ext^2_\Ga(J(\ov k),\ov k^*)&\to &\Br k\\
\end{array} \label{big}
\end{equation}
This diagram is commutative, except for the upper square which is
anti-commutative (with the sign conventions of \cite{Mac}).
The middle and the lower squares are obvious,
so we just need to explain the upper square.
The associativity of the Yoneda product (\cite{Mac}, Ch. III, Thm. 5.3)
implies the commutativity of the upper square if
the maps are the products with the class of (\ref{1}).
By \cite{Mac}, III, Thm. 9.1, such is the left hand vertical map, but
the right hand one differs from the Yoneda product by $-1$.

Let $A$ denote the Albanese variety of $X$.
The abelian varieties $J$ and $A$ are dual to each other.
A choice of a $\ov k$-point on $X$ defines the
Albanese map $\ov X\to\ov A$ over $\ov k$
sending this point to 0. This map canonically
descends to a morphism $X\to D$,
where $D$ is a $k$-torsor of $A$ (cf. \cite{S}, 3.3). Let $\delta(X)\in \H^1(k,A)$
be the class of $D$. This class does not depend on any choice.
In the particular case when $X$ is a $k$-torsor of an abelian
variety, the map $X \to D$ is an isomorphism, so that
$ X(k) \neq \emptyset  $ if and only if  $\delta(X)=0$.

The Barsotti--Weil isomorphism $A(\ov k)=\Ext^{1}_{\ov k-{\rm gps}}(J,\G_m)$
(\cite{GA}, VII.3) gives rise to natural isomorphisms
(\cite{Milne}, Lemma 3.1, p. 50):
\begin{equation}
\H^n(k,A)=\Ext^{n+1}_{k-{\rm gps}}(J,\G_m), \label{BW}
\end{equation}
where $k-{\rm gps}$ is the category of commutative algebraic
groups over $k$, and $n$ is a non-negative integer.
Here the $\Ext ^n$ groups are defined by means of equivalence classes
of  $n$-extensions.

Building upon these isomorphisms, one defines two Tate pairings.

The first Tate pairing
$$\H^1(k,J) \times A(k) \to \Br k$$
is defined by means of the composition of maps
$$A(k) = \Ext^{ 1}_{k-{\rm gps}}(J,\G_m) \to \Ext^1_\Ga(J(\ov k),\ov k^*) \to \Hom(\H^1(k,J), \Br k),$$
where the first map is the isomorphism (\ref{BW})  for $n=0$, the second map is the forgetful map,
the third map is the Yoneda pairing.

The second Tate pairing
$$J(k) \times \H^1(k,A) \to \Br k$$
is defined by means of the composition of maps
$$\H^1(k,A) =  \Ext^{ 2}_{k-{\rm gps}}(J,\G_m) \to \Ext^2_\Ga(J(\ov k),\ov k^*) \to \Hom(J(k), \Br k),$$
where the first map is the isomorphism (\ref{BW}) for $n=1$, the second map is the forgetful map,
the third map is the Yoneda pairing.

A legitimate question, which we need not address, is whether these two pairings coincide upon
swapping    $A$ with $J$.
As the referee points out, biextensions should help.

The second Tate pairing fits into the commutative diagram
\begin{equation}
\begin{array}{ccccc}
(\Pic\ov X)^\Ga&\times&\Ext^2_\Ga(\Pic\ov X,\ov k^*)&\to &\Br k\\
\uparrow&&\downarrow&&||\\
J(k)&\times&\Ext^2_\Ga(J(\ov k),\ov k^*)&\to &\Br k\\
||&&\uparrow&&||\\
J(k)&\times&\H^1(k,A)&\to &\Br k\\
\end{array} \label{v}
\end{equation}
where the  top square comes from the diagram (\ref{big}) (the pairing in the middle being the Yoneda pairing).

\bpr \label{prop:e-delta}
In this diagram,
the image of $e(X) \in \Ext^2_\Ga(\Pic\ov X,\ov k^*)$ in $\Ext^2_\Ga(J(\ov k),\ov k^*)$
is equal to the image of $\delta(X) \in \H^1(k,A)$ in $\Ext^2_\Ga(J(\ov k),\ov k^*)$.
\epr

This is  \cite{S1}, Prop. 2.1.

\subsection{The Brauer group and the elementary obstruction over local fields}  \label{sec:local-general}

Let $R$ be a
henselian, discrete, rank one  valuation ring with finite residue field
and field of fractions $k$ of characteristic zero.
We shall  here refer to such a field as henselian local field
(for $k$ of arbitrary characteristic, see Milne \cite{Milne},
 I.2, p. 43). A henselian local field is a $p$-adic field if and only if
it is complete.

\bthe \label{t}
Let $X$ be a geometrically integral variety
over a henselian local field $k$.
Then $ob(X)=0$ if and only if
the natural map $\Br k\to \Br k(X)$ is injective.
\ethe

{\em Proof\/}
Over any field, the assumption $ob(X)=0$ implies that
$\Br k \to \Br k(X)$ is injective
(Lemma \ref{facts1} (iii)).

Using resolution of singularities
we may assume $X$ smooth and projective.
Assume that $\Br k\to \Br k(X)$ is injective.
This implies that $\Br k \to \Br X$ is injective,
and hence
the map $(\Pic\ov X)^\Ga\to\Br k$ in sequence (\ref{e1}) is zero.
This map is the cup-product with $e(X)$ (Lemma \ref{facts2} (iii)),
thus $e(X)$ is orthogonal to $(\Pic\ov X)^\Ga$
with respect to the Yoneda product.

Consider  the   diagram (\ref{big}).
Now $({\rm Pic\,} \ov X)^\Ga$ is orthogonal to
$e(X) \in  \Ext^2_\Ga({\rm Pic\,} \ov X,\ov k^*)$, thus
the image of $e(X)$ in $\Ext^2_\Ga(J(\ov k),\ov k^*)$
is orthogonal to $J(k)$.
As recalled in Proposition \ref{prop:e-delta}, this image
is equal to the image of $\delta(X)$
under the bottom right hand vertical map in diagram (\ref{v}).
{From} that diagram, we conclude that
$\delta(X) \in \H^1(k,A)$ is orthogonal to $J(k)$
under the second Tate pairing.
 By  Tate's second duality  theorem (\cite{Milne}, I.3, Thm. 3.2 (statement for $\alpha^2$),  Cor. 3.4 and Remark 3.10, l. 5 on p.~59)
this implies
$\delta(X)=0$. Hence the image of
$e(X) \in  \Ext^2_\Ga({\rm Pic\,} \ov X,\ov k^*)$
in $\Ext^2_\Ga(J(\ov k),\ov k^*)$ is zero. Thus
$e(X)$ is the image of some element
$g(X) \in \Ext^2_\Ga(\NS \ov X,\ov k^*)$.
This element is orthogonal to the image of
$({\rm Pic\,} \ov X)^\Ga$ in $(\NS \ov X)^\Ga$.
Let $M \subset \H^1(k,J)$ be the image of $(\NS \ov X)^\Ga$.
Since the abelian group $\NS \ov X$ is finitely generated, and
$\H^1(k,J)$ is torsion, the abelian group $M$ is finite.
The cup-product with $g(X)$ defines a map
$$ (\NS \ov X)^\Ga \to \Br k=\Q/\Z$$
which induces a map
$\nu : M \to \Q/\Z$. Since $\Q/\Z$ is an injective abelian group,
the following natural homomorphism is surjective:
\begin{equation}
\Hom_\Z(\H^1(k,J),\Q/\Z)\to \Hom_\Z(M,\Q/\Z). \label{surj}
\end{equation}

As explained above,
the Barsotti--Weil isomorphism (\ref{BW})
$A(k)= \Ext^1_{k-{\rm gps}}(J,\G_{m})$ and the forgetful map
$\Ext^1_{k-{\rm gps}}(J,\G_{m}) \to  \Ext^1_\Ga(J(\ov k), \ov k^*)$
give rise to the diagram
\begin{equation}
\begin{array}{ccccc}
\H^1(k,J) & \times & A(k) & \to & \Br k \\
 ||&&\downarrow&&||\\
\H^1(k,J)&\times & \Ext^1_\Ga(J(\ov k),\ov k^*)&\to &\Br k,\\
\end{array}\label{Tatepairing}
\end{equation}
which is the definition of the upper row pairing (\cite{Milne}, Prop. 0.16 p.~14 and   I.3):
this is the first Tate pairing as defined at the end of the previous subsection.
By Tate's first duality theorem over a henselian local field
(\cite{Milne}, I.3, Thm. 3.2, statement for $\alpha^1$, Cor. 3.4 and Remark 3.10, l. 5 on p.~59),
this pairing
induces a perfect duality between the discrete group $\H^1(k,J)$
and the completion $A(k)^{\hat{}}$ of $A(k)$
with respect to the natural topology on $k$. In particular,
$A(k)$ is a dense subgroup of $\Hom_\Z(\H^1(k,J),\Q/\Z)$.
By the surjectivity of (\ref{surj}),
its image in $\Hom_\Z(M,\Q/\Z)$ is also dense.
Thus the image of $A(k)$
is the whole  finite set $\Hom_\Z(M,\Q/\Z)$.
Hence there exists an element of $A(k)$ which induces $\nu$ on $M$
via the first Tate pairing.
Let $\rho \in \Ext^1_\Ga(J(\ov k),\ov k^*)$ be its image.
If one modifies $g(X) \in \Ext^2_\Ga(\NS \ov X,\ov k^*)$ by the image of $\rho$ under the map
$\Ext^1_\Ga(J(\ov k),\ov k^*) \to \Ext^2_\Ga(\NS \ov X,\ov k^*)$, one obtains
an element $g_{1}(X) \in \Ext^2_\Ga(\NS \ov X,\ov k^*)$ whose image in
$\Ext^2_\Ga(\Pic \ov X, {\ov k}^*)$ is still $e(X)$, but which is
now orthogonal to $(\NS \ov X)^\Ga$ with respect to the cup-product pairing
$$(\NS \ov X)^\Ga \times \Ext^2_\Ga(\NS \ov X,\ov k^*) \to \Br k.$$
The N\'eron--Severi group $\NS \ov X$ is a discrete Galois module of finite type.
Over the henselian local field $k$, the latter pairing defines
an isomorphism between the groups $\Ext^2_\Ga(\NS \ov X,\ov k^*)$ and
$\Hom_{\Z}((\NS \ov X)^\Ga, \Q/\Z)$ (\cite{Milne}, I.2, Thm. 2.1 and 2.14).
Thus $e(X)=0$. QED.

\bigskip

{\it Remarks\/}
(1) Let $X$ be a smooth, projective, geometrically integral $k$-variety.
Recall that the existence of a zero-cycle of degree 1 on $X$ implies $ob(X)=0$
(Lemma \ref{facts1}  (ii)). If $X$ is a curve over a $p$-adic field,
the converse is also true by a theorem
of Roquette and Lichtenbaum \cite{L}. For $X$ of arbitrary dimension over a $p$-adic field, under the assumption that
$X$ has a regular model ${\cal X}$ proper over the ring of integers of $k$,
one conjectures the equivalence of the two statements:

(a) There exists a zero-cycle of degree 1 on $X$.

(b) The map $\Br k \to \Br X/ \Br {\cal X}$ is injective.

\noindent It is known (\cite{CSai}, Thm. 3.1)   that (a) implies (b), and
that (b) implies the existence of a
zero-cycle of degree a power of $p$. The proof of this last result given
in \cite{CSai} was conditional upon the conjectured absolute purity for the prime-to-$p$ part of the Brauer group of ${\cal X}$; that property is now known, thanks to results of  Gabber (see \cite{Fu}).

(2)
Over a $p$-adic field $k$, for any integer  $n \leq 8$,
there exist  smooth cubic hypersurfaces $X \subset \P^n_k$
 which have no rational point, hence (D. Coray \cite{coray}) no zero-cycles of degree 1.
 If the dimension
of the hypersurface is at least 3, Lemma \ref{facts2} (ix) gives  $ob(X)=0$.

(3) The  theorem as it stands does not extend to arbitrary fields $k$ of cohomological dimension 2.  Let $k=\Cc(u,v)$ be the rational function field in 2 variables.
The quadric $Q\subset\P^3_{k}$ given by
 $$ X^2+uY^2+vZ^2+(1+u)uvT^2=0$$
has no $k$-points, as one sees by going over to $\Cc((u))((v))$,
but it satisfies $\Br k \hookrightarrow \Br k(Q)$. For
$K=k({\sqrt{1+u}})$ the group $\Br K$ does not inject into $\Br K(Q)$.
For more on this example see Subsection \ref{7}

\medskip

Recall that a field $R$ is real closed if $-1$ is not a sum of squares in
$R$, but is a sum of squares in any finite extension of $R$.
By the Artin--Schreier theorem
$[\ov R:R]=2$.

\bthe \label{treal}
Let $X$ be a geometrically integral variety
over a real closed field $R$.
Then $ob(X)=0$ if and only if
the natural map $\Br R \to \Br R(X)$ is injective.
\ethe

{\em Proof\/} The proof is the same as the proof given above, once
one takes into account the following two results.

Let $A$ and $B$ be dual abelian varieties over the field $R$.
Let $C$ be the algebraic closure of $R$.
The natural pairing
$$A(R) \times \H^1(R,B) \to \Br R = \Z/2 \subset \Q/\Z$$
induces a perfect pairing of finite $2$-torsion groups
$$ A(R)/N_{C/R}A(C) \times \H^1(R,B) \to \Q/\Z$$
(over $R=\R$ see \cite{Milne}, I.3, Remark 3.7; in the general case,
see \cite{Geyer}).

Let $\Ga=\Gal(C/R)$. Let $M$ be a finitely generated $\Ga$-module.
Then the natural pairing
$$M^\Ga \times \Ext^2_\Ga(M, C^*) \to \Br R = \Z/2$$
induces an isomorphism
$$\Ext^2_\Ga(M, C^*) \simeq \Hom_{\Z}(M^\Ga/N_{C/R}M, \Z/2)$$
(see \cite{Milne}, I.2, Thm. 2.13; the proof is given for $R=\R$
but it holds for an arbitrary real closed field). QED

\medskip

{\it Remark\/} It is easy to give examples of varieties $X$
over an arbitrary real closed field $R$ such that
$ob(X)=0$ but $X(R)=\emptyset$, for example anisotropic quadrics in $\P^n$
for $n \geq 4$. It is however known that a smooth, geometrically integral $R$-variety
$X$ has an $R$-point if and only if for all $i$ the maps $\H^i(R,\Z/2) \to \H^i(R(X),\Z/2)$
are injective. That the first statement implies the second is a general fact for
smooth varieties over a field, with a rational point, which may be seen in a number of ways.
If $X/R$ is geometrically integral of dimension $d$  and has no $R$-point, then the cohomological dimension of the field $R(X)$ is equal to $d$. This is a consequence of a theorem of Serre
(see \cite{CTP}, Prop. 1.2.1). For modern developments of this classical topic, see \cite{Scheiderer}.

\medskip

We give a short, new proof of the following theorems of
J. van Hamel (\cite{vHR}, Section 5 for $k$ the field of real numbers and
\cite{vH} for  $k$ a $p$-adic  field). This theorem
 generalizes   previous results  of Roquette and of Lichtenbaum.

\bthe[van Hamel] \label{le1}
Let $X$ be a smooth, projective, geometrically integral variety
over a henselian local field $k$
or over a real closed field.
Then $ob(X)=0$ implies
$\delta(X)=0$.
In particular, a $k$-torsor $X$ of an abelian
variety is trivial if and only if $ob(X)=0$.
\ethe

{\em Proof\/}
Consider the diagram (\ref{v}). As recalled in
Proposition \ref{prop:e-delta},
the image of $e(X) \in \Ext^2_\Ga(\Pic\ov X,\ov k^*) $ in
$\Ext^2_\Ga(J(\ov k),\ov k^*)$
is equal to the image of $\delta(X) \in \H^1(k,A)$ in
$\Ext^2_\Ga(J(\ov k),\ov k^*)$.
The hypothesis $ob(X)=0$ implies $e(X)=0$ (Lemma \ref{facts2} (i)).
Hence $J(k)$ is orthogonal to $\delta(X)$ with respect
to the bottom pairing of (\ref{v}). Since $k$ is a henselian local field
or a real closed field,
Tate's second duality theorem implies that $\delta(X)=0$. QED

\medskip

Let $k$ be a henselian local field and let ${\hat k}$ be its completion.
The following lemma is well known.

\ble \label{S0}
Let the fields $k$ and $\hat k$ be as above. The natural map $\Br k\to \Br{\hat k}$ is an isomorphism.
\ele

The following result goes back to  \cite{Gr}.
\bpr \label{greenberg} Let the fields $k$ and $\hat k$ be as above.  If a $k$-algebra of finite type
admits a $k$-algebra homomorphism to $\hat k$, then
it admits a $k$-algebra homomorphism to $k$. In particular,
the field $\hat k$ is the union of its $k$-subalgebras of
finite type $A$ admitting a retraction $A \to k$.
\epr

This implies that for any contravariant functor $F$  from $k$-schemes to sets
which commutes with filtering  limits  with affine transition morphisms
the natural map $F(X) \to
F(X\times_k{\hat k})$ is injective. This applies in particular to the functor
  $F(X)=\Br X$. This  also implies that for any $k$-variety $X$
   the  conditions $X(k)\neq \emptyset$
  and $X(\hat k) \neq \emptyset$ are equivalent.

\bpr \label{H}
Let the fields $k$ and $\hat k$ be as above.
Let $X$ be a smooth geometrically integral variety over
$k$. Then $ob(X)=0$ if and only if $ob(X\times_k {\hat k})=0$.
\epr

{\em Proof\/} The previous comment implies that the map
$\Br X \to \Br (X\times_k{\hat k})$ is injective.
Together with Lemma \ref{S0}, this shows that $\Br k \to \Br X$
is injective if and only if $\Br{\hat k} \to
\Br (X\times_k {\hat k})$ is injective. In turn, this implies that
$\Br k \to \Br k(X)$ is injective if and only if $\Br {\hat k}
\to \Br {\hat k}(X)$ is injective. A double application of
Theorem \ref{t} completes the proof. QED

\medskip

Let now $k \subset  R$ be an inclusion of real closed fields.
The analogue of Greenberg's result is a classical theorem going back
to E. Artin: if a $k$-algebra of finite type admits a $k$-homomorphism
to $R$, then it admits a $k$-homomorphism to $k$.
The natural map $\Br k \to \Br R =\Z/2$ is a bijection.
Theorem \ref{treal} and the same argument as above now give:

\bpr \label{Hreal}
Let  $k \subset  R$ be an inclusion of real closed fields.
Let $X$ be a smooth geometrically integral variety over
$k$. Then $ob(X)=0$ if and only if $ob(X\times_k {R})=0$.
\epr

\subsection{The Brauer group and the elementary obstruction over number fields}\label{sec:global-general}

\bpr \label{obob}
Let $X$ be a smooth geometrically integral variety over
a number field $k$, and $k_v$ be the completion of $k$ at a  place $v$.
Then $ob(X)=0$ implies $ob(X\times_k k_v)=0$.
\epr

{\em Proof\/} Let $\tilde k$ be the integral closure of $k$ in $k_{v}$.
For $v$ finite, this is the fraction field of the
henselization of the ring of integers of $k$ at $v$.
For $v$ real,
this is a real closed field.
Since $\tilde k\subset\ov k$, the condition $ob(X)=0$
implies $ob(X\times_k\tilde k)=0$. Now the statement follows
from Propositions \ref{H} and \ref{Hreal}. QED
\medskip

Recall that by definition
$$
{\cyr B}(X)=\Ker [\Bro  X\to \prod_v \Bro  X_v/\textup{Br}_0 X_v],
$$
where $\textup{Br}_0 X_v$ is the image of $\Br k_v$ in $\Bro X_v$.
This group does not change under restriction of $X$ to a nonempty open set
(\cite{San}, Lemme 6.1).

\bigskip

Recall that $X(\A_{k})^{\cyr B}$
is the subset of $X(\A_{k})$ consisting of the adelic points
orthogonal to ${\cyr B}(X)$ with respect to the Brauer--Manin
pairing (see \cite{S}, 5.2, for more details). Obviously,
this set  either is empty  or coincides with $X(\A_{k})$.

\bthe \label{B}
Let $X$ be a smooth, geometrically integral variety
over a number field $k$.
Assume that $X(\A_{k}) \neq \emptyset$ and
$ob(X)=0$. Then $X(\A_{k})=X(\A_{k})^{\cyr B}$. In particular,
$X(\A_{k})^{\cyr B} \neq \emptyset$.
\ethe

{\em Proof\/} Let us fix a Galois-equivariant section $\sigma$ of the map
$\k^* \to \k(X)^*$.
For each place $v$ of $k$ fix a decomposition group $\Ga_{v } \subset
\Ga ={\rm Gal}(\k/k)$.
Let ${\tilde k}_{v} \subset \k$ be the fixed field of $\Ga_{v}$.
If $v$ is finite, this is a henselian local field.
 If $v$ is
a real place of $k$, then this is a real closure of $k$.
Let $\alpha \in {\cyr B}(X)$.
For each place $v$ of $k$, the image of $\alpha$ in $\Br X_{v}$
comes from a well defined element of $\Br k_{v} $.
Using the same arguments as in the end of the previous subsection,
we see that the restriction of $\alpha$ to
$\Br (X\times_k{{\tilde k}_{v}})$
  comes from a well defined element $\xi_{v}$ of $\Br {\tilde k}_{v}$.
  This last element may be computed  by composing the maps
    $$  \Bro(X\times_k{{\tilde k}_{v}})  \to \H^2(\Ga_{v},\k(X)^*)   \to
\H^2(\Ga_{v},\k^*),$$
    where the last map is given by $\sigma$.
    We also have the element $\xi \in \Br k$ which is the image of $\alpha$
    under the composite map $\Bro X \to \H^2(\Ga,\k(X)^*) \to \H^2(\Ga,\k^*)$,
    where the last map is induced by $\sigma$. Now $\xi_{v}$ is
clearly the restriction
    of $\xi \in \Br k$ to $\Br {\tilde k}_{v}$.
Thus the sum of the local invariants
associated to
    the family $\xi_{v}$ is the sum of the local invariants of
$\xi$, it is therefore
zero. QED.

\medskip

{\it Remark\/} We keep the assumption $X(\A_{k}) \neq \emptyset$.
In the particular case when $\Pic\ov X$
is a free abelian group, a delicate theorem asserts that the conditions $ob(X)=0$ and
$X(\A_{k})^{\cyr B}=X(\A_{k})$
are equivalent (\cite{CS}, Prop. 3.3.2).
It would be interesting to see if the same is true in general\footnote{Wittenberg \cite{Witt}, building upon
work of Harari and Szamuely \cite{HS}, has now proved: if one grants the finiteness of Tate--Shafarevich groups
of abelian varieties over number fields, then the answer to this question is positive.}.

\bigskip
\def\Pico{\textup{Pic}^0}

We conclude this subsection with the following observation, which  does not seem to be documented in literature (but see \cite{Mazur},  Cor. 1, p.~40, for a similar result).
\bpr
Let $X$ be a smooth, proper, geometrically integral variety over a number field $k$, and  let $A=\Pico{X}$ be its Picard variety. Assume that for any finite extension $K/k$
the Tate--Shafarevich group of $A_{K}$ is finite. Then
the quotient of ${\cyr B}(X)$ by the image of $\Br k$ is finite.
\epr

{\em Proof\/}
 We have the exact sequence of Galois modules
 $$0 \to \Pico {\ov X} \to \Pic {\ov X}  \to {\rm NS}\, {\ov X}  \to 0.$$
 Let  $K/k$ be a finite Galois extension such that $X(K) \neq \emptyset$ and  the composite map
  $$\Pic X_{K} \to \Pic {\ov X}  \to  {\rm NS}\, {\ov X} $$ is onto.
   Let
  $\mathfrak{h}$ be the Galois group of  $\k$ over $K$.  The  $\mathfrak{h}$-module  ${\rm NS}\, {\ov X} $
 is the direct sum of a free abelian group
  $\Z^r$ and a finite abelian group $F$, both with trivial action of $\mathfrak{h}$.
Galois cohomology yields the exact sequence
$$0 \to \H^1(K, \Pico {\ov X}  )  \to \H^1(K,  \Pic {\ov X} )  \to \H^1(K,F).$$
We have analogous exact sequences over each henselization ${\tilde K}_w$ of $K$:
 $$0 \to \H^1({\tilde K}_w,\Pico  {\ov X} ) \to \H^1({\tilde K}_w, \Pic  {\ov X} ) \to \H^1({\tilde K}_w,F). $$
 By Chebotarev's theorem, the
 kernel of the diagonal map
$\H^1(K,F) \to \prod_{w} \H^1({\tilde K}_w,F)$, where $w$ runs through all places   of $K$, vanishes.
By our assumption on Tate--Shafarevich groups, the kernel of
$ \H^1(K,\Pico {\ov X}  ) \to \prod_{w} \H^1({\tilde K}_w,\Pico {\ov X} )$
is finite. Thus the kernel of
$\H^1(K, \Pic {\ov X} ) \to \prod_{w}\H^1({\tilde K}_w, \Pic {\ov X} )$ is finite.

Let $G$ be the finite group ${\rm Gal}(K/k)$. We have the
standard restriction-inflation exact sequence
$$0 \to \H^1(G, (\Pic {\ov X} )^\mathfrak{h}) \to \H^1(k,\Pic {\ov X} ) \to \H^1(K,\Pic {\ov X} ).$$
The Mordell--Weil theorem and the N\'eron--Severi theorem imply that
the abelian group  $\Pic X_{K} = (\Pic {\ov X} )^\mathfrak{h}$ is of finite type.
Thus  $\H^1(G,(\Pic {\ov X} )^\mathfrak{h}) $ is finite.
It is then clear that the kernel of
$\H^1(k, \Pic {\ov X} ) \to \prod_{v} \H^1({\tilde k}_{v}, \Pic {\ov X} )$ is finite.

The  argument given in the proof of Theorem  \ref{B} shows that the group
${\cyr B}(X)$ may also be defined as
$$
{\cyr B}(X)=\Ker [\Bro  X\to \prod_v \Bro  X_{\tilde{k}_{v}} /\textup{Br}_0 X_{\tilde{k}_{v}}],
$$
where $\textup{Br}_0 X_{\tilde{k}_{v}}$ is the image of $\Br {\tilde k}_v$ in $\Bro X_{\tilde{k}_{v}}$.

 From the Hochschild--Serre spectral sequence for the multiplicative group and the projection map
$X \to {\rm Spec}\, k$ we have the standard exact sequences
$$0 \to \textup{Br}_0 X \to \textup{Br}_{1} X \to \H^1(k, \Pic {\ov X} )$$
and for each place $v$ of $k$
$$0 \to \textup{Br}_0 X_{\tilde{k}_{v}} \to \textup{Br}_{1} X_{\tilde{k}_{v}} \to \H^1(\tilde{k}_{v}, \Pic {\ov X} ).$$

The group ${\cyr B}(X)/ \textup{Br}_0 X$ is thus a subgroup of the kernel of the diagonal map
$\H^1(k, \Pic {\ov X} ) \to \prod_{v} \H^1({\tilde k}_{v}, \Pic {\ov X} )$. It is thus finite. QED

\section{Homogeneous spaces}

By convention, all homogeneous spaces we shall consider will be right homogeneous spaces.

\subsection{Structure of algebraic groups}\label{sec:structure}

Let $k$ be a field of characteristic 0.

The following theorem will be constantly used. If $H \hookrightarrow G$ is a homomorphism of
(not necessarily affine) algebraic groups over $k$ which is an immersion, then the quotient
$G/H$ exists in the category of $k$-varieties (A. Grothendieck, \cite{Groth},  Thm. 7.2, Cor. 7.4;
P.~Gabriel, \cite{Ga}, Thm. 3.2 p.~302).

We shall  also use the  fact: if $H \subset G$ is a normal subgroup of an algebraic group over $k$, and $X$ is a $k$-variety which is a right homogeneous space of $G$, then the quotient variety
$Y=X/H$ exists  in the category of $k$-varieties, it is a (right) $G/H$-homogeneous space.
The morphism $X \to Y$ is faithfully flat and smooth.
When $G$ is affine, a proof of this fact is given
in \cite{B96}, Lemma 3.1.  By the above result of  Grothendieck and Gabriel,
that proof works for arbitrary algebraic groups.

If $L$ is a connected linear group, we denote by
$L\uu$ its unipotent radical, a normal connected subgroup of $L$.
We let $L\red$ be the quotient of $L$ by its unipotent radical $L\uu$.
This is a connected reductive group.
We let $L\sss \subset L\red$ be the derived group of $L\red$. This is
a connected semisimple group.
We denote by
$L\tor$ the biggest toric quotient of $L$.
The kernel of  $L \to L\tor$ is a normal, connected subgroup of $L$ denoted by $L\ssu$.
The group $L\ssu$ is an extension of  $ L\sss$ by $L\uu$.

Any connected algebraic group $G$ over $k$ is an extension
\begin{equation}
1 \to L \to G \to A \to 1 \label{LGA}
\end{equation}
of an abelian variety $A/k$ by a normal, connected linear $k$-group $G$
(Chevalley's theorem \cite{Ros, Conrad}).
We write $L=G\lin$.
This is a characteristic subgroup of $G$, it is stable under all automorphisms of the group $G$.
We denote   by  $Z(G)$   the centre of $G$  and
by $G\sab$  the biggest group quotient of
$G$ which is a semiabelian variety.
We write $G^\der$ for the derived subgroup $[G,G]$.
The group $G^\der$ is clearly contained in $L$, hence is a connected linear
algebraic group.

 If $L$ is reductive, then $L^\der=G^\der$ hence in particular $G^\der$ is a semisimple group.
   Indeed,
the connected semisimple  group $L^\der$ is normal in $G$, the quotient  $G'$ of $G$ by
$L^\der$ is an extension of $A$ by the group $L/L^\der$, which
is  $L\tor$.
Any  group extension of an abelian variety by a torus is central.
Since there are no nonconstant morphisms from an abelian variety to a torus,
any  such  group extension is commutative.
Thus $G'$  is a semiabelian variety.
Since the kernel  $L^\der$ of $G \to G'$ is semisimple, we have $G'=G\sab$
and $L^\der=G^\der$.

By Prop. 4 of  \cite{Wu} the connected group
$G/Z(G)$ is linear. According to \cite{Wu}, Thm.~1, we have
the following commutative diagram:
\begin{equation}
\begin{array}{ccccccccc}
1&\to &Z(L)&\to &Z(G)&\to &A&\to &0\\
&&\downarrow&&\downarrow&&||&&\\
1&\to &L&\to &G&\to &A&\to &0
\end{array}
\label{e4}
\end{equation}

Let $H$ be a linear $k$-group (not necessarily connected).
We write $\hat H$ for the group of characters of $\ov H$ (this is a finitely generated
discrete Galois module), and
$H^{\rm mult}$ for the biggest quotient of $H$ which is
a $k$-group of multiplicative type. By construction
the $k$-groups $H$ and $H^{\rm mult}$ have the same groups of characters.
We set $$H_1=\ker[H\to H^{\rm mult}].$$
In Theorems \ref{thm:non-connected},   \ref{thm:main-n}  and
 \ref{thm:non-connected-Manin}
we shall
make the hypothesis that
$H _1$ is connected and that $\hat H_1=0$.
This hypothesis is satisfied if $H$ is connected.
Indeed, in this case the group $H_{1}$ coincides with the connected group $H\ssu$,
and clearly $\Hbar\ssu$ has no nontrivial characters.
For general $H$ the hypothesis need not be satisfied:
consider the example where $H$ is a finite,
noncommutative, solvable group, or the case of a noncommutative extension
of a finite abelian group by a torus.

\bpr \label{p1}
Let $X$ be a homogeneous space of a connected $k$-group
whose maximal connected linear subgroup
has trivial unipotent radical. Assume that the stabilizers
of the geometric points of $X$ are connected.
Then $X$ can be given the structure of a homogeneous space of an algebraic group $G$
satisfying the following conditions:

$G^{\rm lin}$ has trivial unipotent radical,

$G^\der$ is semisimple simply connected,

the stabilizers of
the geometric points of $X$ in $G$ are linear and connected.
\epr

{\em Proof\/}
Let $G$ be a connected group whose maximal linear subgroup $L$
has trivial unipotent radical. Assume that $G$
acts transitively on $X$ with connected geometric stabilizers.
The group $G$ is an extension (\ref{LGA}).
According to (\ref{e4}) we have $L/Z(L)=G/Z(G)$. Since $L$
is reductive, the latter group is semisimple. This also implies
$G^\der=L^\der$, as explained above.

We write $\St_{\ov x, \ov G}$ for the stabilizer of
$\ov x\in X(\ov k)$ in $\ov G$. These subgroups of $\ov G$
form one conjugacy class.
\medskip
 
{\it First reduction.}

The subgroup $Z(\ov G)\cap \St_{\ov x, \ov G}$ is central in $\ov G$,
and does not depend on $\ov x$. Hence it is stable under
the action of the absolute Galois group
$\Ga$, and so
$Z(\ov G)\cap \St_{\ov x, \ov G}=\ov C$ for
a central subgroup $C\subset G$. Then
$X$ is a homogeneous space of $G/C$ such that
$\St_{\ov x,\ov G/\ov C}=\St_{\ov x, \ov G}/\ov C$.
The group $G/Z(G)$ is linear, hence so is $\St_{\ov x, \ov G}/\ov C$.
Replacing $G$ by $G/C$ we may thus assume
 without loss of generality
 that the stabilizers of the
geometric points are linear and connected.
\medskip

{\it Second reduction.}

It is well known  (Langlands, see \cite{MS}, Prop. 3.1) that
given the connected reductive group $L/Z(L)$ there exist
exact sequences of connected reductive
algebraic groups
$$ 1 \to S \to H \to L/Z(L) \to 1$$
with $S$ a $k$-torus central in $H$, and
$H^\der$ simply connected.
(Such extensions are called $z$-extensions.)
 Define $G'$ as the fibred product of $G$ and $H$
over $L/Z(L)$, so that there is a commutative diagram of
exact sequences of algebraic groups
$$
\begin{array}{ccccccccc}
&&&&1&&1&&\\
&&&&\downarrow&&\downarrow&&\\
&&&&S&=&S&&\\
&&&&\downarrow&&\downarrow&&\\
1&\to &Z(G)&\to &G'&\to &H&\to &1\\
&&||&&\downarrow&&\downarrow&&\\
1&\to &Z(G)&\to &G&\to &L/Z(L)&\to &1\\
&&&&\downarrow&&\downarrow&&\\
&&&&1&&1&&
\end{array}
$$
Note that $Z(G)$ is in the centre of $G'$.
We then have the commutative diagram of exact sequences
of connected linear algebraic groups
$$
\begin{array}{ccccccccc}
&&1&&1&&&&\\
&&\downarrow&&\downarrow&&&&\\
&&S&=&S&&&&\\
&&\downarrow&&\downarrow&&&&\\
1&\to &L'&\to &G'&\to &A&\to &0\\
&&\downarrow&&\downarrow&&||&&\\
1&\to &L&\to &G&\to &A&\to &0\\
&&\downarrow&&\downarrow&&&&\\
&&1&&1&&&&
\end{array}
$$
where $L'$ is the kernel of the composite map $G' \to G \to A$.
Clearly $L'$ is linear, so it is the maximal linear subgroup of $G'$.
Thus the natural map $L'/Z(L')\to G'/Z(G')$ is an isomorphism
of semisimple groups.
Since $Z(G)$ is a central subgroup of $G'$,
the map $G'\to G'/Z(G')$ factors as $G' \to H \to G'/Z(G')$.
The maps $L' \to G' \to H \to G'/Z(G')$
give rise to a series of maps
\begin{equation}
(L')^\der \to (G')^\der \to H^\der \to (G'/Z(G'))^\der=L'/Z(L'),\label{s}
\end{equation}
where the composite map is induced
by the natural map $L' \to L'/Z(L')$.
Since $L'$ is reductive, the first map in (\ref{s}) is an isomorphism,
as was explained above. The maps $G' \to H \to G'/Z(G')$
are surjective, hence so are the second and the third maps in (\ref{s}).
Since $L'$ is a reductive group, the natural map $(L')^\der \to L'/Z(L')$
is an isogeny, hence $(L')^\der \to H^\der$ is also an isogeny.
But $H^\der$ is simply connected since
$H$ is a $z$-extension. This forces $(L')^\der \simeq H^\der$,
so that $(L')^\der=(G')^\der$
is a semisimple simply connected group.
Replacing $G$ by $G'$ we keep the
property that the stabilizers of the geometric points are
connected linear groups. QED

\subsection{Local fields: semiabelian varieties}

\bthe \label{le2}
Let $k$ be a henselian local field or a real closed field.
A $k$-torsor $X$ of a semiabelian variety is trivial
if and only if $ob(X)=0$.
\ethe

{\em Proof\/} Let $X$ be a torsor of
a semiabelian variety $G$,
an extension of an abelian variety $A$ by a torus $T$:
$$1\to T\to G\to A\to 0.$$
Let $D$ be the quotient of $X$ by the action of $T$;
this is a $k$-torsor of $A$, which can also be defined as the
push-forward of $X$ with respect to the map $G\to A$.
By functoriality $ob(D)=0$, so that $D\simeq A$
by Theorem \ref{le1}. Thus $X$ is an $A$-torsor of $T$.
We write $\xi$ for the class of this torsor in $\H^1_{\et}(A,T)$,
and $\xi_m\in \H^1(k,T)$ for the class of the fibre $X_m$
at a $k$-point $m$ of $A$. Our goal is to find $m$ with $\xi_m=0$.

From the bilinear pairing of $k$-group schemes
$$\hat T  \times   T  \to \G_{m,k}$$
we deduce a cup-product pairing
$$\H^1(k,{\hat T}) \times \H^1_{\et}(A,T) \to \H^2_{\et}(A,\G_{m})=\Br A.$$
the image of which lies in $\Br_{1} A$.

Let $B\subset\Bro   A$ be the subgroup consisting
of the elements $\alpha\cup\xi$, where $\alpha\in \H^1(k,\hat T)$.
The group $\H^1(k,\hat T)$ is finite, hence so is $B$.

The $k$-point $0\in A(k)$ defines a splitting of (\ref{e1}) applied to $X=A$,
so that $\Bro   A$ decomposes as the direct sum of $\Br k$
and the subgroup consisting of the elements
${\cal A}\in\Bro   A$ such that ${\cal A}(0)=0$,
naturally identified with $\H^1(k,\Pic\ov A)$.
The canonical map
$r : \Bro   A\to\H^1(k,\Pic\ov A)$ can be written as
${\cal A}\mapsto {\cal A}-{\cal A}(0)$.
Let $J$ be the Picard variety of $A$, which is also
the dual abelian variety of $A$.

We now prove the following statements, the last of which
proves the theorem:
\medskip

(1) The restriction of the canonical map
$r : \Bro   A\to\H^1(k,\Pic\ov A)$
to $B$ factors through $\H^1(k,J)$.

(2) $B\cap\Br k=0$.

(3) There exists $m\in A(k)$ orthogonal to $B$ with respect to the
pairing $A(k)\times\Br A\to\Br k$ given by the evaluation.

(4) For any point $m$ satisfying (3) we have $\xi_m=0$, that is,
the fibre of $X\to A$ over $m$ contains a $k$-point.

\medskip

\noindent {\em Proof of\/} (1). Let $\lambda : \hat T\to \Pic\ov A$
be the type of the torsor $X\to A$ (\cite{CS},  (2.0.2);  \cite{S}).
It is well known (see \cite{GA}, chap. VII, no. 16, Thm. 6 and  comment thereafter) that $\ov X$ can be given
the structure of a group extension of $\ov A$ by $\ov T$ if and only if
$\lambda$ factors through the natural injection
$J(\ov k)\hookrightarrow \Pic\ov A$. Now (1) follows from
Thm. 4.1.1 of \cite{S}, which says that
the following diagram commutes
\begin{equation}
\begin{array}{ccc}
\H^1(k,\hat T) &\hfl{\lambda_*}{}{8mm} &\H^1(k,J)\\
\vfl{\cup\xi}{}{4mm}&&\vfl{}{}{4mm}\\
\Bro   A&\hfl{r}{}{8mm} &\H^1(k,\Pic\ov A)
\end{array}
\label{e3}
\end{equation}

\noindent {\em Proof of\/} (2).
The image of $\xi$ under the base change map
$\H^1(A,T)\to\H^1(X,T)$ is zero since $X\times_A X$
is a trivial $X$-torsor (the diagonal is a section).
Thus $B$ goes to 0 under the pull-back map $\Br A\to\Br X$.
The assumption $ob(X)=0$ implies that the natural map $\Br k\to\Br X$
is injective, and this implies $B\cap\Br k=0$.

\medskip

\noindent {\em Proof of\/} (3).
We now {\it define} a pairing
\begin{equation}
A(k)\times \H^1(k,\hat T)\to \Br k, \label{nor}
\end{equation}
in the following manner. A couple $(m,\alpha) \in A(k)\times \H^1(k,\hat T)$
is sent to
$$(\alpha\cup\xi)(m)-(\alpha\cup\xi)(0)=\alpha\cup(\xi_m-\xi_0).$$
 We claim that this pairing is bilinear.
To prove this, consider the diagram of pairings
\begin{equation}
\begin{array}{ccccc}
A(k)&\times&\H^1(k,J)&\to&\Br k\\
||&&\downarrow&&||\\
A(k)&\times&\H^1(k,\Pic\ov A)&\to&\Br k
\end{array} \label{b}
\end{equation}
where the top row is the
Tate
 pairing,
and the bottom row is the pairing given by
evaluating elements of $\H^1(k,\Pic\ov A)$,
understood as the subgroup of $\Bro   A$ consisting of
the elements with trivial value at $0$.
This diagram commutes by Prop. 8(c) of \cite{M}.
{From} (\ref{e3}) we see that the map
$\H^1(k,\hat T)\to \H^1(k,\Pic\ov A)$ sending $\alpha$
to $r(\alpha\cup \xi)=\alpha\cup\xi-(\alpha\cup\xi)(0)$ factors through
$\H^1(k,J)$, and so the pairing (\ref{nor})
is bilinear since such is the top pairing
of (\ref{b}).

We now use the hypothesis on the field $k$.
There is a natural embedding $\Br k \hookrightarrow \Q/\Z$.
The pairing (\ref{nor}) induces a homomorphism
$\sigma : A(k)\to B^*=\Hom(B,\Q/\Z)$.
Let us show that $\sigma$ is surjective. If it is not,
there exists $b\in B$, $b\not=0$, such that
$\sigma(m)$ applied to $b=\xi \cup \alpha \in \Bro A$ is zero for any $m$,
that is, $b(m)-b(0)=0$ for all $m\in A(k)$. Thus
 $b-b(0) $
comes from an element of $\H^1(k,J)$ orthogonal
to $A(k)$ with respect to the Tate pairing.
However,
over a henselian local field or   a real closed field
the right kernel of the Tate pairing $A(k) \times \H^1(k,J) \to \Br k$
is
zero, hence $b=b(0) \in B \subset \Bro A$ is a non-zero constant element
in $B$. This contradicts (2).

By the surjectivity of $\sigma$ there exists $m\in A(k)$
such that $\sigma(m)$ is the element of $B^*$ given by
$b\mapsto -b(0)$ for any $b\in B$. This says that $b(m)-b(0)=-b(0)$,
so that $b(m)=0$ for any $b\in B$. This finishes the proof of (3).

\medskip

\noindent {\em Proof of\/} (4). By (3) we
have $(\alpha\cup\xi)(m)=\alpha\cup\xi_m=0$ for any
$\alpha\in \H^1(k,\hat T)$. Hence $\xi_m$
is orthogonal to $\H^1(k,\hat T)$
with respect to the pairing
$$\H^1(k,\hat T)\times \H^1(k,T)\to\Br k.$$
For  $k$
a henselian local field or a real closed field,
 this pairing is non-degenerate (\cite{Milne}, I.2, Thm. 2.14 (c) and Thm. 2.13, whose proof works over
 a real closed field)
 thus $\xi_m=0$.
This finishes the proof of the theorem. QED

\subsection{$p$-adic   fields: main theorem}\label{5}

\bthe \label{t2}
Let $X/k$ be a homogeneous space of a connected $k$-group
(not necessarily linear) such that the stabilizer $\Hbar$
of a geometric point  $\xbar\in X(\kbar)$ is connected.
If $k$ is a henselian local field,
 then $X$ has a $k$-point
if and only if $ob(X)=0$.
\ethe

In conjunction with Theorem \ref{t} this gives the following corollary.

\begin{cor} \label{cor2}
Let $X/k$ be a homogeneous space of a connected $k$-group
(not necessarily linear) such that the stabilizer $\Hbar$
of a geometric point   $\xbar\in X(\kbar)$ is connected.
If $k$ is a henselian local field,
then $X$ has a $k$-point
if and only if $\Br k $ injects into $\Br k(X)$.
\end{cor}
\medskip

{\it Proof of  Theorem  $\ref{t2}$}.
{\it First reduction}

Suppose that $X$ is a right homogeneous space of a
connected group $G$ represented as an extension (\ref{LGA}).
The unipotent radical $L\uu\subset L$ is a normal subgroup
of $G$. Let
$G'=G/L\uu$. This group satisfies $(L')\uu=0$.
The following properties, proved in Lemma 3.1 of \cite{B96},
hold over any perfect field $k$.
The quotient $X'=X/L\uu$
exists, and there is a natural
projection map $X \to X'$. This map is surjective on $\ov k$-points
and its geometric fibres are orbits of $L\uu$.
The variety $X'$ is a homogeneous space of $G'$ with connected
geometric  stabilizers.

The hypothesis $ob(X)=0$ implies $ob(X')=0$.
Suppose we have found a $k$-point $y \in X'(k)$.
Then the fibre $X_y$ is a $k$-variety which is
a homogeneous space of the unipotent $k$-group $L\uu$.
According to Lemma 3.2(i) of \cite{B96}, over any perfect field $k$ this
implies
$X_y(k)\neq \emptyset$. Thus $X(k) \neq \emptyset$.

Thus without loss of generality we may assume that
the unipotent radical $L\uu$ of $L$ is trivial, so that $L$ is reductive.

\medskip

{\it Second reduction}

By Proposition \ref{p1} we can further assume that $G^\der$
is semisimple simply connected, and the stabilizers of
the geometric points of $X$ in $G$ are linear and connected.
This reduction has nothing to do with the nature of the field $k$.
It does not change $X$, hence we keep the assumption $ob(X)=0$.

\medskip
{\it Relaxing the assumptions}

To prove Theorem  \ref{t2} it is  enough to prove
the following result (whose proof is similar to that of Thm. 2.2 in \cite{B96}).
We write $G\sss$ for $L\sss$, and $G\uu$ for $L\uu$, where $L=G\lin$.
The notation $\Hbar_1$ was defined in Subsection \ref{sec:structure}

\bthe \label{thm:non-connected}
Let $k$ be a henselian local field,
$G$ a connected $k$-group, and $X/k$ a
homogeneous space of $G$ with geometric stabilizer $\Hbar$. Assume

{\rm (i)} $G\uu=\{1\}$,

{\rm (ii)} $\Hbar \subset \ov G\lin$,

{\rm (iii)} $G\sss$ is simply connected,

{\rm (iv)} $\Hbar_1$ is connected and
has no non-trivial characters (e.g. $\Hbar$ is connected).

\noindent Then $ob(X)=0$ if and only if $X(k) \neq \emptyset$.
\ethe

The homogeneous space $X$ defines a $k$-form of $\Hbar^{\rm mult}$
which we denote by $M$ (see \cite{B96}, 4.1).
We have a canonical homomorphism $M\to G\sab$.
For this, see the computation at the end of  Subsection 1.2 of \cite{B99}.
In that paper, $G=L$ is linear, the
  calculation   uses  the commutativity of $L\tor$.
It generalizes to the present context  with the commutative group $G\sab$ in place of $L\tor$.

Here is another way to construct
the homomorphism $M\to G\sab$. One extends the base field from $k$ to the function field $k(X)$
of $X$. Consideration of the stabilizer $H'$ of the generic point of $X$
yields a map $H' \to  G\times_kk(X)$ over $k(X)$ which induces a map
$H' \to G\sab\times_k{k(X)}$. Since  $H$ hence $H'$ is linear, this map factors through
$T\times_kk(X)$, where $T=(G\sab)\lin$ is the maximal torus inside the semiabelian variety
 $G\sab$.
There is then an induced   $k(X)$-morphism  $M\times_kk(X) \to T\times_kk(X)$. Such a map
comes from a unique morphism $M \to T \subset G\sab$.

We first prove a special case of Theorem \ref{thm:non-connected}.

\bpr \label{5.4}
With the hypotheses of Theorem $\ref{thm:non-connected}$,
assume that $M$ {\bf injects} into $G\sab$
(i.e. $\Hbar\cap {\ov G}\sss=\Hbar_1$).
Then $X$ has a $k$-point.
\epr

{\em Proof\/}
Set $Y=X/G\sss$.
Then $Y$ is a homogeneous space of the semiabelian variety $G\sa$,
hence it is a torsor of some semiabelian variety.
We have a canonical map $X\to Y$.
{From}  $ob(X)=0$ we deduce $ob(Y)$=0 (see the beginning of the introduction).
 By Theorem \ref{le2},
 $Y$ has a $k$-point $y$.
 Let $X_y$ denote the fibre of $X$ over $y$.
It is a homogeneous space of $G\sss$
with geometric stabilizer $\Hbar\cap {\ov G}\sss=\Hbar_1$.
The group $G\sss$ is semisimple simply connected by (iii).
The group $\Hbar_1$ is connected and has no nontrivial characters by (iv).
By \cite{B93}, Thm. 7.2
(that theorem is stated over a $p$-adic field, but it also   holds over a henselian local field,
see  the proof of Theorem \ref{good} hereafter)
the $k$-variety $X_y$ has a $k$-point.
Hence $X$ has a $k$-point. QED
\medskip

For the general case we need an easy lemma.

\ble\label{lem:killH2}
Let $M$ be a $k$-group of multiplicative type
and $\eta\in \H^2(k,M)$ a cohomology class.
Then there exists an embedding $j\colon M\into P$
into a quasi-trivial $k$-torus $P$
such that $j_*(\eta)=0$.
\ele

{\em Proof\/}
We can embed $M$ into a quasi-trivial torus, and so
assume without loss of generality
that $M=R_{K/k}\Gm$ for some finite extension $K/k$.
We have a canonical isomorphism $s_K\colon \H^2(k,R_{K/k}\Gm)\isoto \H^2(K,\Gm)$.
Let $L/K$ be a finite extension such that the image of $s_K(\eta)$
in $\H^2(L,\Gm)$ is zero.
Consider the natural injection of quasi-trivial tori
$c_{K/L}\colon R_{K/k}\Gm\into R_{L/k}\Gm$.
Then $(c_{K/L})_*(\eta)=0$. QED

\medskip

Let us resume the proof of Theorem \ref{thm:non-connected}.
Let $\xbar\in X(\kbar)$ be a point with stabilizer $\Hbar$.
Let $\eta_X\in \H^2(k,\Hbar,\kappa)$ be the cohomology class defined by $X$
(Springer's class, see \cite{B93}, 7.7, or \cite{Sp}, 1.20),
where $\kappa$ is the $k$-kernel defined by $X$, see \cite{B93}, 7.1.
Recall that $\Hbar_1=\ker[\Hbar\to\Hbar^\mult]$.
Clearly the subgroup $\Hbar_1$ is invariant under all semialgebraic automorphisms of $\Hbar$,
hence $\kappa$ induces a $k$-kernel $\kappa^\mult$ in $\Hbar^\mult$,
and we obtain a map
$$
\mu_*\colon \H^2(k,\Hbar,\kappa)\to \H^2(k, \Hbar^\mult,\kappa^\mult)
$$
induced by the canonical map $\mu\colon\Hbar\to\Hbar^\mult$, see \cite{B93}, 1.7.
Since $\Hbar^\mult$ is an abelian group, $\kappa^\mult$ defines a $k$-form of $\Hbar^\mult$,
which is the $k$-form $M$ mentioned above.
We obtain an element $\mu_*(\eta_X)\in \H^2(k,M)=\H^2(k,\Hbar^\mult,\kappa^\mult)$.
Note that in \cite{B93}, section 7,   $G$ is assumed semisimple and simply connected,
but the general constructions we refer to hold with $G$ any $k$-group,
the key point is that the subgroup $\Hbar$ is linear.

By Lemma \ref{lem:killH2}
we can construct an embedding $j\colon M\into P$
into a quasi-trivial $k$-torus $P$ such that $j_*(\mu_*(\eta_X))=0$.
Consider the $k$-group $F=G\times P$, and the embedding
$$
\Hbar\into \ov F=F\times_k\ov k\quad \text{given by} \quad h\mapsto (h, j(\mu(h))).
$$
Set $\Zbar=\Hbar\backslash \ov F$.
We have a right action $\abar\colon \Zbar\times \ov F\to \Zbar$
and an $\ov F$-equivariant map
$$
\pibar\colon \Zbar\to \ov X,\quad \Hbar\cdot (g,p)\mapsto \Hbar\cdot g,
\text{\quad where }g\in \Gbar,\ p\in \overline{P}.
$$
Then $\Zbar$ is a homogeneous space of $\ov F$ with respect to the action $\abar$,
and the map $\pibar\colon \Zbar\to \ov X$ is a torsor under ${\ov P}$.
The homomorphism $M\to F\sa$ is injective.

In \cite{B96}, 4.7, it is proved that
$\text{Aut}_{\ov F,\ov X}(\Zbar)=P(\kbar)$.
By \cite{B96}, Lemma 4.8, the element $j_*(\mu_*(\eta_X))\in \H^2(k, P)$
is the only obstruction to the existence
of a $k$-form $(Z,a,\pi)$ of the triple $(\Zbar, \abar,\pibar)$:
there exists such a $k$-form if and only if $j_*(\mu_*(\eta_X))=0$.
In our case by construction we have $j_*(\mu_*(\eta_X))=0$, hence
there exists a $k$-form $(Z,a,\pi)$ of  $(\Zbar, \abar,\pibar)$.
Since $\pi\colon Z\to X$ is a torsor under the quasi-trivial torus $P$,
from Hilbert's theorem 90 and Shapiro's lemma
we conclude that $Z$ is $k$-birationally isomorphic to $X\times P$.
 From  $ob(X)=0$ we deduce $ob(Z)=0$ (see  the beginning of the introduction).

We obtain a homogeneous space $Z$ of a connected $k$-group $F$
such that  $F\sss$ is simply connected, with geometric stabilizer $\Hbar$.
The group $M$ injects into $F\sa=G\sa\times P$, and $ob(Z)=0$.
By Proposition \ref{5.4},
 $Z$ has a $k$-point.
Thus $X$ has a $k$-point. QED.

\bigskip

{\it Remark\/} In Thm. 3.9 of \cite{Fl}, M. Florence constructs a homogeneous space
$X$ of $G=\PGL(D)$, for a  quaternion algebra $D$ over a $p$-adic field, such that
the geometric stabilizer ${\overline H}\simeq\Z/2 \times \Z/2$, and
$X$ has a zero-cycle of degree 1
but no rational points. The space $X$ can also be viewed
as a homogeneous space of $\SL(D)$,
the geometric stabilizer now being the quaternion group.
Since $X$ has a zero-cycle of degree 1,
the map $\Br k \to \Br k(X)$ is injective. Thus $ob(X)=0$. This shows that in Theorem
\ref{thm:non-connected} neither condition (iii) nor condition (iv)
may
be omitted.

\subsection{Good fields of cohomological dimension at most $2$} \label{7}

A field of characteristic zero is called a good field of cohomological dimension
at most $2$ if it satisfies the following properties:

(i) Its cohomological dimension $cd(k)$ is at most 2.

(ii) Over any finite field extension $K/k$, for any central simple algebra $A/K$,
the index of $A$ (as a $K$-algebra) and the exponent of the class of $A$ in $\Br K$ coincide.

(iii) For any semisimple simply connected group $G/k$
we have $\H^1(k,G)=0$.

\medskip

According to Serre's ``Conjecture II'', (i) should imply (iii). This is known for groups of classical type.
The combination of (i) and (ii) implies (iii) for all groups without
factors of type $E_{8}$ (see the references in \cite{CGP}).

Properties (i) to (iii) are satisfied for
henselian local fields
and
for totally imaginary number fields.

For the fraction field of a 2-dimensional strictly henselian local domain, with residue field
of characteristic zero, these three properties also hold  \cite{COP}, \cite{CGP}.

For the function field of an algebraic surface over an algebraically closed field of
characteristic zero, properties (i) and (ii) hold. For (ii), this is de Jong's theorem \cite{dJ}.
Hence in this case  (iii) is known when $G$ has no factors of type $E_{8}$.

\bthe \label{good}
Let $k$ be a good field of cohomological dimension at most $2$
and characteristic zero.
Let $X/k$ be a homogeneous space of a connected {\bf linear} group $G$.
Assume that the geometric stabilizers are connected.
Then $X(k) \neq \emptyset$ if and only if $ob(X)=0$.
\ethe

{\em Proof\/}  We follow the proof of Theorem \ref{t2}.
The first and second reduction have nothing to do with the nature of the field $k$.
It remains to prove the analogue of Theorem \ref{thm:non-connected}.
Since $G$ here is linear,
the semiabelian variety $G\sab$ is a $k$-torus. With the notation as in
the proof of Proposition \ref{5.4},
the $k$-variety $Y$ is a homogeneous space of a $k$-torus. It satisfies $ob(Y)=0$.
Over any field, this implies $Y(k) \neq \emptyset$, see Lemma \ref{facts1} (iv).
Keeping the notation of Proposition \ref{5.4},
we find $y \in Y(k)$, and then the $k$-variety $X_{y}$
is a homogeneous space of $G\sss$
with geometric stabilizer $\Hbar\cap {\ov G}\sss=\Hbar_1$.
The group $G\sss$ is semisimple simply connected.
The group $\Hbar_1$ is connected and has no nontrivial characters.
   Over a good field of cohomological dimension 2,
the analogue of \cite{B93}, Thm. 7.2,
 is  Propositions 5.3 and 5.4 of \cite{CGP},
which build upon  the key Theorem 2.1 of \cite{CGP}
and use the formalism of \cite{B93}.
This shows that the $k$-variety $X_y$ has a $k$-point.
Hence $X$ has a $k$-point. This completes the proof of the analogue of Proposition \ref{5.4}.

Lemma \ref{lem:killH2} holds over any field.
The rest of the proof of Theorem \ref{thm:non-connected}
is a reduction to Proposition \ref{5.4}, which works equally well over
any ground field. QED

\bco
Let $k$  be a good field of cohomological dimension at most $2$
and characteristic zero.
Let $X/k$ be a homogeneous space of a connected linear group $G$.
Assume that the geometric stabilizers are connected.

{\rm (i)} Then $X(k) \neq \emptyset$  if and only if for any flasque $k$-torus $S$,
the restriction map $\H^2(k,S) \to \H^2(k(X),S)$ is injective.

{\rm (ii)} If $X$ is projective, then $X(k) \neq \emptyset$ if and only if for any finite field extension $K/k$
the map $\Br K \to \Br K(X)$ is an injection.

{\rm (iii)}
If $X$ is projective and the abelian group $\Pic(\ov X)$ is free of rank $1$,
then $X(k) \neq \emptyset$ if and only if
the natural map $\Br k \to \Br k(X)$ is an injection.
\eco

 {\em Proof\/}
 (i) This follows from   Theorem A of \cite{CK}, Theorem \ref{good},
 and \ref{facts2} (iv).

(ii) The Bruhat decomposition implies that
the geometric Picard group of a projective homogeneous space
of a connected linear group is a permutation $\Ga$-module
(cf. \cite{CGP}, the proof of Lemma 5.6 on p. 337).
Now (ii) follows from Theorem \ref{good} and Lemma \ref{facts2} (vi).

(iii) This follows from Theorem \ref{good} and Lemma \ref{facts2} (v).
QED

 \medskip

Remark 3 after Theorem \ref{t}
  shows that in (ii)
above one cannot simply assume the injectivity of
$\Br k \to  \Br k(X)$.

 \bigskip

 {\it Remarks}
 (1) For any even integer $n=2m \geq 6$,
Merkurjev \cite{Me} constructs a (big) field $k_{n} $ of cohomological dimension 2
and an anisotropic quadratic form of rank $n$ over $k_{n}$.
The associated quadric is a homogeneous space of a spinor group with connected
geometric stabilizers. There are elements of order 2 in the
Brauer group of $k_n$ which are not of index 2.
Thus the mere hypothesis $cd(k) \leq 2$ is not enough for the above theorem
 to hold,
condition (ii)
(in the definition of a good field of cohomological dimension  at most $2$)
 is required.

   \medskip

 (2) The above corollary should be compared with
the recent work of de Jong and Starr \cite{dJSt}
on projective homogeneous varieties over function fields in two variables.

  \medskip

(3) Let $k=\Cc(u,v)$ be the rational function field in two variables over the complex field.
Let $X \subset \P^8_{k}$ be the smooth cubic hypersurface defined by the diagonal cubic
form with coefficients $1,u,u^2,v,vu,vu^2,v^2,v^2u,v^2u^2$.
One easily checks that $X(k)=\emptyset$. In fact, $X$ has no points in
$\Cc((u))((v))$.
On the other hand, Lemma 2.2 (ix)
ensures $ob(X)=0$.
The same comment applies to smooth cubic hypersurfaces in $\P^n_{k}$ with $4 \leq n \leq 7$
  defined by taking subforms of the above form.

\subsection{Number fields}

Let $k$ be a number field.
We write $\Omega_{r}$ for the set of all {\it real}  places of $k$.
We set $k_\infty=\prod_{v\in \Omega_{r}} k_v$, then for a $k$-variety $X$
we have $X(k_\infty)=\prod_{v\in\Omega_{r}} X(k_v)$.
When $k$ is totally imaginary, the following result is a special case of
Theorem \ref{good}.

\bthe \label{real}
Let $k$ be a number field and $X/k$ be a homogeneous space of a connected {\bf linear} algebraic $k$-group $G$
with connected geometric stabilizer. Assume that $X$ has a $k_v$-point for every real place
$v$ of $k$.
If $ob(X)=0$, then $X$ has a $k$-point.
\ethe

Proceeding as in Subsection \ref{5}, we see that this is a consequence of the following result,
whose proof is similar to that of Theorem 2.2
in \cite{B96}.

\bthe\label{thm:main-n}
Let $k$ be a number field and $X/k$ be a homogeneous space of a connected {\bf linear} algebraic $k$-group $G$
with geometric stabilizer $\Hbar$.
Assume that:
\par {\rm (i)} $G\uu=\{1\}$,
\par {\rm (ii)} $G\sss$ is simply connected,
\par {\rm (iii)} $\Hbar_1$ is connected and has
no non-trivial characters,
\par {\rm (iv)} $X$ has a $k_v$-point for every $v\in\Omega_{r}$.
\par\noindent If $ob(X)=0$, then $X$ has a $k$-point.
\ethe

The homogeneous space $X$ defines a $k$-form of $\Hbar^{\mult}$
which we denote by $M$.
We have a canonical homomorphism $M\to G\tor$.
We first prove a special case of Theorem \ref{thm:main-n}.

\bpr\label{prop:injects-n}
In Theorem $\ref{thm:main-n}$
assume that $M$ injects in $G\tor$
(i.e. $\Hbar\cap {\ov G}\sss=\Hbar_1$).
Then $X$ has a $k$-point.
\epr

{\em Proof\/}
Set $Y=X/G\sss$.
Then $Y$ is a homogeneous space of the $k$-torus $G\tor$,
hence it is a torsor of some $k$-torus $T$.
We have a canonical map $\alpha\colon X\to Y$.
Since $ob(X)=0$, we see that $ob(Y)$=0.
Hence $Y$ has a $k$-point $y$ by Lemma \ref{facts1} (iv).

  The map $\alpha\colon X\to Y$ is smooth, hence for $v\in \Omega_{r}$
the image $\sY_v:=\alpha(X(k_v))$ is open in $Y(k_v)$ and nonempty (because $X$ has a $k_v$-point).
Set $\sY_\infty=\prod_{v\in\Omega_{r}} \sY_v$, then $\sY_\infty$ is a nonempty open subset in $Y(k_r)$.
By the real approximation theorem for tori (due to J-P. Serre), see \cite{San}, Cor. 3.5, or
\cite{Vos}, Thm. 11.5, the set $Y(k)$ is dense in $Y(k_\infty)$.
Hence there exists a $k$-point $y'\in Y(k)\cap \sY_\infty$.

Consider the fibre $X_{y'}$ of $X$ over $y'$.
It is a homogeneous space of $G\sss$
with geometric stabilizer $\Hbar\cap {\ov G}\sss =\Hbar_1$.
The group $G\sss$ is semisimple simply connected by (ii).
The group $\Hbar_1$ is connected and has no nontrivial characters by (iii).
Since $y'\in\sY_\infty$, the variety $X_{y'}$ has a $k_v$-point for every  $v\in\Omega_{r}$.
By \cite{B93}, Thm. 7.3 (vi) and Cor. 7.4,  $X_{y'}$ has a $k$-point.
Hence $X$ has a $k$-point. QED
\medskip

We resume the proof of Theorem \ref{thm:main-n}.
Let $G$ and $X$ be as in that theorem.
Let $\xbar\in X(\kbar)$ be a point with stabilizer $\Hbar$.
We have a canonical map
$\mu_*\colon \H^2(k,\Hbar,\kappa)\to \H^2(k,M)$, where $\kappa$ is the $k$-kernel
defined by $X$.
Let $\eta_X\in \H^2(k,\Hbar,\kappa)$ be the cohomology class defined by $X$.
Consider $\mu_*(\eta_X)\in \H^2(k,M)$.
By Lemma \ref{lem:killH2}
we can construct an embedding $j\colon M\into P$
into a quasi-trivial $k$-torus $P$ such that $j_*(\mu_*(\eta_X))=0$.

As in the proof of Theorem \ref{thm:non-connected},
we construct the $k$-group
$F=G\times P$, and a triple $(Z,a,\pi)$, where $(Z,a)$ is a homogeneous space of $F$
and $(Z,\pi)$ is a torsor of $P$ over $X$.
Since $(Z,\pi)$ is a torsor of the quasi-trivial torus $P$ over $X$,
and $X$ has a $k_v$-point for any $v\in \Omega_{r}$,
we see that $Z$ has a $k_v$-point for such $v$.
Also since $(Z,\pi)$ is a torsor of the quasi-trivial torus $P$,
we see that $Z$ is $k$-birationally isomorphic to $X\times P$.
Since $ob(X)=0$, we see that $ob(Z)=0$.

We obtain a homogeneous space $Z$ of a connected reductive $k$-group $F$
such that  $F\sss$ is simply connected, with geometric stabilizer $\Hbar$.
The group $M$ injects into $F\tor=G\tor\times P$, and $ob(Z)=0$.
The homogeneous space  $Z$ has a $k_v$-point for any $v\in \Omega_{r}$.
By Proposition \ref{prop:injects-n}, $Z$ has a $k$-point.
Thus $X$ has a $k$-point. QED

\medskip

{\it Remark} To prove Theorem \ref{real},
one could also argue as follows.
According to Proposition \ref{obob},
the hypothesis $ob(X)=0$
implies $ob(X\times_{k}k_{v})=0$ for each nonarchimedean place $v$ of $k$.
Theorem \ref{t2}
then implies $X(k_{v}) \neq \emptyset$ for each
nonarchimedean
 place $v$.
Thus $X(\A_{k}) \neq \emptyset$. Theorem \ref{B}
 then implies
 $X^{c}(\A_{k})^{\cyr B}=X^{c}(\A_{k}) \neq \emptyset$.
 From Theorem 2.2 of \cite{B96}
  we conclude $X(k) \neq \emptyset$.
 This proof looks more elegant than the one above,
but it relies on Theorem 2.2 of \cite{B96}, whose proof occupies most of the paper \cite{B96}.
In the proof given above, one sees precisely where  the   linearity of   $G$   is used. It is  to ensure
weak approximation at the real places for $Y$, which is   a principal homogeneous space of a
torus (a similar argument occurs in  \cite{B96}). Had we not assumed $G$ linear, $Y$ would have been a principal homogeneous space of
a semiabelian variety. For an abelian variety $A$ over a number field $k$, weak approximation at the real places may  badly  fail: over some real completion $k_{v}$, there may be no $k$-point in a  connected component of $A(k_{v})=A(\R)$. This will be the basis of the example given in  Subsection 3.6.

\medskip

The question as to whether the Brauer--Manin obstruction attached to
${\cyr B}(X)$ is the only obstruction to the Hasse principle
on $k$-torsors of arbitrary connected algebraic groups
was raised in \cite{S} (p.~133, Question 1).
D. Harari and T. Szamuely  \cite{HS} recently announced a positive solution
to this problem for torsors of
semiabelian varieties.

\bthe[Harari--Szamuely] \label{hs}
Let $k$ be a number field, and $X$ a $k$-torsor
of a semiabelian variety $G$. Assume that the Tate--Shafarevich
group of the biggest quotient of $G$ which is an abelian variety,
is finite. If $X$ has a family of local points $P_v\in X(k_v)$,
for all places $v$ of $k$, which is
orthogonal to ${\cyr B}(X)$ with respect to the Brauer--Manin pairing,
then $X$ has a $k$-point.
\ethe

This implies the following global analogue of Theorem \ref{t2}.

\bthe \label{imago}
Let $X$ be a homogeneous space of a (not necessarily linear) connected group $G$
such that the stabilizers
of the geometric points of $X$ are connected.
Assume that the Tate--Shafarevich
group of the biggest quotient of $G$ which is an abelian variety,
is finite. If
$k$ is a totally imaginary number field, then $X$ has a $k$-point
if and only if $ob(X)=0$.
\ethe
{\em Proof\/} We follow the proof of Theorem \ref{t2}
up to the place where Theorem \ref{le2} is used, and
apply Theorems \ref{B} and \ref{hs} instead.
Thm. 7.2 (local) and Cor. 7.4 (global) of
\cite{B93}
allow us to finish the proof in the same way as before. QED

\subsection{Number fields: an example}

We now proceed to construct a $\Q$-torsor $X$
of a non-commutative connected  algebraic group over $\Q$, such that $ob(X)=0$,
$X$ has points over all completions of $\Q$,  further $X^{c}(\A_{\Q})^{\cyr B}=X^{c}(\A_{\Q}) \neq \emptyset$,  but   $X$ has no
$\Q$-points. Thus in general the answer to the aforementioned question
is negative.

\medskip

Let $E/\Q$ be the elliptic curve with affine equation
$$y^2=(x^2-3)(x-2).$$
We note that the set $E(\R)$ has two connected components:
the connected component of the origin of the group law,
given by $x\geq 2$, and the component given by $x^2\leq 3$.

The quaternion algebra $(x-2,-1)$ over $\Q(E)$ comes from a (unique) Azumaya algebra
over $E$, which will be denoted by $A$. If $M$ is a $p$-adic or a real
point of $E$, then the value of $A$ at $M$ is either $0$, or
the unique element of $\Br\Q_v$ of order 2.

An application
of {\tt magma} shows that
$E(\Q)=\{0,(2,0)\}$, but in what follows we shall only need
the following statement.

\ble \label{la}
For any prime $p$ and any point $M_p\in E(\Q_p)$
the value $A(M_p)$ is zero.
The sum $\sum_v A(M_v)$, taken over all places $v$ of $\Q$,
is zero if and only if $M_\R$ is in
the connected component of $0 \in E(\R)$. In particular,
$E(\Q)$ is contained in the connected component of $0 \in E(\R)$.
\ele
{\em Proof\/} We first prove that $A$ takes only trivial values on
$\Q_p$-points of $E$, for any prime $p$.
It is enough to compute the values of $A$ at the points
$M_p=(x,y)$ such that $xy\neq 0$. Indeed since $A$
is an Azumaya algebra over $E$, for each place $v$ of $\Q$,
the map $E(\Q_{v}) \to \Z/2$ given by evaluation of $A$ is continuous,
and for any nonempty Zariski open set $U$ of $E$,
$U(\Q_{v})$ is dense in $E(\Q_{v})$.
Let $K=\Q(\sqrt{-1})$.

Let $p$ be an odd prime.
If $p$ splits in $K$, i.e. if $p\equiv 1$ mod 4, then
$-1$ is a square in $\Q_{p}$ and the assertion is trivial.
If  $p$ is inert in $K$, i.e. $p\equiv 3$ mod 4, then
$\alpha\in\Q_{p}^*$ is a norm from $K_{p}$, which is
equivalent to
$(\alpha,-1)=0 \in \Br\Q_{p}$, if and only if $v_{p}(\alpha)$ is even.
If $v_{p}(x)<0$, then
$2v_{p}(y)=v_{p}((x^2-3)(x-2))=3v_{p}(x)$.
Hence $v_{p}(x)$ is even, and then $v_{p}(x-2)$ is even, and so
$(x-2,-1)=0 \in \Br\Q_{p}$.
Assume $v_{p}(x-2)\geq 0$. If $v_{p}(x-2)>0$, then $v_{p}(x^2-3)=0$. Hence
$2v_{p}(y)=v_{p}(x-2)$, so that $v_{p}(x-2)$ is even, and we conclude as before.

Let $p=2$.
Write $x=u/v$ with $u \in \Z_{2}$ and $v \in \Z_{2}$, not both divisible by 2.
In $\Z_{2}$ we have a relation
\begin{equation}
z^2= (u^2-3v^2)(uv-2v^2)\not=0. \label{*}
\end{equation}
If $(u,v)\equiv (0,1)$ or $(1,0)$ mod 2,
then $u^2-3v^2\equiv 1$ mod 4.
In both cases we find $(u^2-3v^2,-1)=0 \in \Br\Q_{2}$.
{From} (\ref{*}) we conclude that $(x-2,-1)=0 \in \Br\Q_{2}$.
It remains to consider the case  $(u,v)\equiv (1,1)$ mod 2.
Write $x=1+2n$ with $n \in \Z_{2}$. Then $x-2=-1+2n$
and $x^2-3= -2+4n+4n^2$. Thus $(x-2)(x^2-3)= 2+4m$ for some
$m \in \Z_{2}$ and this cannot be a square. So there are no such points $(x,y)$.

Finally, if $(x,y) \in E(\R)$, $y\not=0$,
then $(x-2,-1)_{\R}=0$ is equivalent to $x>2$.
Using reciprocity we obtain the statement about
$E(\Q)$. QED
\medskip

Let $f : E'\to E$ be the unramified double covering given by
$u^2=x-2$. The curve $E'$ has a $\Q$-point above $0$;
choosing it for the origin of the group law on $E'$ turns $f$
into an isogeny of degree 2.
We note that $f(E'(\R))$ is the connected component of $0$
of $E(\R)$, so that $E(\R)/f(E'(\R))=\Z/2$.

Let $D$ be the Hamilton quaternions. The group
$L=\SL_1(D)$ is a $\Q(\sqrt{-1})/\Q$-form of $\SL_2$, in
particular, it is semisimple simply connected
with centre $\{\pm 1\}$. Define $G=(\SL_1(D)\times E')/(\Z/2)$,
where $\Z/2$ is generated by $(-1,P)$, $P\in E'(\Q)$,
$f(P)=0$, $P\not=0$.
We obtain a commutative diagram of
extensions of algebraic groups
\begin{equation}
\begin{array}{ccccccccc}
1&\to &\Z/2&\to &E'&\to &E&\to &0\\
&&\downarrow&&\downarrow&&||&&\\
1&\to &\SL_1(D)&\to &G&\to &E&\to &0
\end{array}
\label{e}
\end{equation}
This gives rise to the following commutative diagram
of pointed sets
\begin{equation}
\begin{array}{ccccc}
E(\Q)&\to &\H^1(\Q,\Z/2)&\to &\H^1(\Q,E')\\
||&&\downarrow&&\downarrow\\
E(\Q)&\to &\H^1(\Q,\SL_1(D))&\to &\H^1(\Q,G)
\end{array}
\label{es}
\end{equation}
and the  compatible diagrams with $\Q_{p}$ or $\R$ in place of $\Q$:

We have the canonical isomorphisms
$$\H^1(\Q,\Z/2)=\Q^*/\Q^{*2}, \quad \H^1(\R,\Z/2)=\R^*/\R_{>0},$$
$$\H^1(\Q_{p},\SL_1(D))=\Q_p^*/{\rm Nrd}((D\otimes_{\Q}\Q_{p})^*)=1$$
$$\H^1(\Q,\SL_1(D))=\H^1(\R,\SL_1(D))=\R^*/\R_{>0}.$$
The map $\Z/2 \to \SL_1(D)$ induces a surjection $\Q^*/\Q^{*2} \to \H^1(\Q,SL_1(D))$
which itself  induces a bijection $\{\pm 1\} = \H^1(\Q,SL_1(D))$.

In the above diagrams, the map  $E(\Q) \to \Q^*/\Q^{*2}$
  on the affine open set of $E$ defined by
$x-2 \neq 0$ is given by evaluation of the function $x-2$.
As one easily checks, the value on the point at infinity is $1$, the value on
the point $x=2$  is  the  value taken  by $x^2-3$, namely $1$.
The same statement holds over any field extension of $\Q$.

 \bpr \label{example}
Let $G/\Q$ be the above defined algebraic group. Let $X$
be a torsor of $G$ whose class $\xi \in \H^1(\Q,G)$ is the image of  $-1 \in \H^1(\Q,\Z/2)$
under the map $$\H^1(\Q,\Z/2) \to
\H^1(\Q,G).$$

Then $ob(X)=0$,
 $X(\A_{\Q})^{\cyr B}=X(\A_{\Q}) \neq \emptyset$
 but $X(\Q)=\emptyset$.

 Let $X^{c}$ be a smooth compactification of $X$.
One has $X^{c}(\A_{\Q})^{\cyr B}=X^{c}(\A_{\Q}) \neq \emptyset$  and
$X^{c}(\A_{\Q})^{\Bro X^{c}}= \emptyset$.
\epr
{\em Proof\/}
We use the commutativity and functoriality of the above diagrams.
From  $\H^1(\Q_{p},\SL_1(D))=1$ we deduce that
 the class of $\xi \in \H^1(\Q,G)$
has trivial image in   $\H^1(\Q_{p},G)$.
From the fact that $E(\R) \to \H^1(\R,\Z/2)$ is onto we deduce that
 the class of $\xi \in \H^1(\Q,G)$
has trivial image $\H^1(\R,G)$.
Thus $X(\A_{\Q}) \neq \emptyset$.

Next, assume that the image of the class of $-1 \in \H^1(\Q,SL_1(D))$
in $\H^1(\Q,G)$ is trivial. Then the image of that
class in $\H^1(\Q,\SL_1(D))$ comes from $E(\Q)$.
Restricting to the cohomology over  $\R$ we see
that the class of $-1$ in $\H^1(\R,\SL_1(D))$,
which is non-trivial, comes from the image
of $E(\Q)$ in $E(\R)$. But $E(\Q)\subset f(E'(\R))$ (Lemma \ref{la})
so this is not possible. Thus $X$ is a non-trivial torsor of $G$,
so that $X(\Q)=\emptyset$.

 Given the torsor $X$ over $\Q$ under the group $G$ we may consider the quotient $Y=X/SL_{1}(D)$
 of $X$ under the action of $SL_{1}(D) \subset G$. This is a torsor over $\Q$ under $E$,
 whose class in $\H^1(\Q,E)$ is the image of $\xi$ under $\H^1(\Q,G) \to \H^1(\Q,E)$.
 The projection map $X \to Y$ makes $X$ into a torsor under $SL_{1}(D)$.
 Since $\xi$ comes from $\H^1(\Q,\Z/2)$ the above diagram shows that the class of $Y$
 in $\H^1(\Q,E)$ is trivial. We may thus identify $Y=E$. All in all, we see that
 $X$ is a torsor over $E$ under $SL_{1}(D)$.

This argument shows that an   open set of $X$ is isomorphic to
the affine variety given by the system of equations
$$ y^2=(x^2-3)(x-2)\neq 0, \quad  2-x=u^2+v^2+w^2+t^2.$$

Let $\Ga=\Gal(\ov \Q/\Q)$. The projection map $X \to E$ induces a Galois equivariant  map from
the   2-extension of
continuous discrete $\Ga$-modules
$$
1\to \ov \Q[E]^*\to \ov \Q(E)^*\to \Div\ov E\to\Pic\ov E\to 0
$$
to the 2-extension
$$
1\to \ov \Q[X]^*\to \ov \Q(X)^*\to \Div\ov X\to\Pic\ov X\to 0.
$$
Over $\ov \Q$, the projection $\ov X  \to \ov E$ makes
$\ov X $ into an $SL_{2}$-torsor over $ \ov E$.
Any such torsor is locally trivial for the Zariski topology.
Any invertible function on $SL_{2}$ is constant, and the Picard group
of the simply connected group $SL_{2}$ is trivial. From this we deduce that
the maps $\ov \Q^* \to \ov \Q[E]^* \to \ov \Q[X]^*$ and $\Pic \ov E \to \Pic \ov X$
are isomorphisms.  Pull-back from $\ov E$ to $\ov X$ thus maps
the 2-extension
$$
1\to \ov \Q^*\to \ov \Q(E)^*\to \Div\ov E\to\Pic\ov E\to 0
$$
to the 2-extension
$$
1\to \ov \Q^*\to \ov \Q(X)^*\to \Div\ov X\to\Pic\ov X\to 0,
$$
the map $\Pic \ov E \to \Pic \ov X$ being an isomorphism.
We have $E(\Q)\neq \emptyset$, hence $ob(E)=0$. Thus  the class of the first extension is trivial,
hence so is
the class of the second extension. This show $ob(X)=0$.

We now have $ob(X^{c})=ob(X)=0$. Theorem \ref{B}
then implies $X (\A_{\Q})^{\cyr B}=X (\A_{\Q}) \neq \emptyset$.
It also implies $X^{c}(\A_{\Q})^{\cyr B}=X^{c}(\A_{\Q}) \neq \emptyset$.
This finishes the proof of the proposition.
QED

\bigskip

{\it Remark}
The computation in Lemma \ref{la} shows that the counterexample
to the Hasse principle on $X$ is due to the Brauer--Manin obstruction
given by $\pi^*A\in\Br    X$.  The class $A \in \Br X$ comes from
$\Br E=\Bro E$, hence lies in
$\Bro X^{c}$. Hence $X^{c}(\A_{\Q})^{\Bro X^{c}}= \emptyset$.
This is  in accordance with a result of Harari
which we shall extend in the Appendix.

\renewcommand{\sV}{{\mathcal{V}}}
\newcommand{\sU}{{\mathcal{U}}}
\renewcommand{\ab}{{}^{\textup{ab}}}

\renewcommand{\thesection} {}
\renewcommand{\thetheorem} {A.\arabic{theorem}}

%\appendix
\section{Appendix: The Brauer--Manin obstruction for homogeneous spaces}

Let $k$ be a number field.
We denote by $\Omega$ the set of all places of $k$,
and by $\Omega_r$ the set of all real places of $k$.
If $S\subset\Omega$, we set $k_S=\prod_{v\in S} k_v$.
If $X$ a $k$-variety, we have $X(k_S)=\prod_{v\in S} X(k_v)$.
In particular, $X(k_\Omega)=\prod_{v\in\Omega}X(k_v)$.

For a connected $k$-group $G$ we write $G\ab:=G/G\lin$,
it is the biggest quotient of $G$ which is an abelian variety.

\begin{theorem} \label{thm:Harari-generalized}
Let $G$ be a connected algebraic group over a number field $k$. Let
$X$ be a  homogeneous space of  $G$ such that the stabilizers of the geometric points of $X$
are connected. Let $X^c$ be a smooth compactification of $X$.
Assume that a point $x_\Omega=(x_v)_{v\in\Omega} \in  X(k_\Omega)$ is orthogonal to
$\Bro X^c$ with respect to the Brauer--Manin pairing.
Assume that the Tate--Shafarevich group of the maximal abelian  variety quotient $G\ab$ of $G$
is finite.  Then for any finite set $S$
of  {\bf nonarchimedean} places of $k$ and any open neighbourhood $\sU_{S}$ of
$x_S=(x_v)_{v \in S}$ in $X(k_S)$ there exists a rational point $x_0 \in X(k)$
whose diagonal image in
$\prod_{v \in S}X(k_{v})$ lies in $\sU_{S}$. Moreover we can ensure  that
for each archimedean place $v$, the points $x_0$ and $x_{v}$ lie in the same connected
component of $X(k_{v})$.
\end{theorem}

This theorem generalizes a recent result of Harari  (\cite{H},
Theorem 1.1), who considers torsors under a connected algebraic
group $G$.
In the extreme case when $G$ is an abelian variety,
our result is due to Manin \cite{M} and Wang \cite{Wang}.
In the other extreme case when $G$ is a linear group,
this result (including approximation at archimedean places)
was obtained in \cite{B96}, Cor. 2.5.
In the general case, a proof by simple devissage
in order to reduce the assertion to these two extreme cases does not work.
Our method of proof uses the reductions and constructions
of Subsections 3.1 and 3.3 in order to reduce the assertion
to the case when $X$ is a $k$-torsor under a semiabelian variety
(treated by Harari \cite{H})
and to the Hasse principle and weak approximation
for  a homogeneous space  of a simply connected semisimple group
with   connected, characterfree  geometric stabilizers (results obtained in  \cite{B93} and \cite{B90}, see also \cite{CTX}).

\medskip

 The proof of Theorem \ref{thm:Harari-generalized}  will occupy the entire appendix.

Let $X$ be a smooth geometrically integral $k$-variety
over a number field $k$.
The Brauer--Manin pairing
$$
X(k_\Omega)\times\Bro X^c\to\Qq/\Zz.
$$
 defines a map
$$
m_X\colon X(k_\Omega)\to(\Bro X^c)^D
$$
where $(\Bro X^c)^D=\Hom(\Bro X^c,\Qq/\Zz)$.
By the birational invariance of the Brauer group \cite{Groth2}, this map does not depend
on the choice of the smooth compactification $X^c$.
 If $\varphi\colon X\to Y$ is a morphism of smooth geometrically integral $k$-varieties,
then by Hironaka's theorem one can construct   smooth compactifications $Y^c$ of $Y$
and $X^c$ of $X$
such that $\varphi$ extends to a morphism $\varphi^c\colon X^c\to Y^c$.
The following diagram  then commutes:
$$
\begin{CD}
X(k_\Omega)  @>{m_X}>> (\Bro X^c)^D \\
@V{\varphi}VV    @VV{\varphi_*}V\\
Y(k_\Omega)  @>{m_Y}>> (\Bro Y^c)^D
\end{CD}
$$
In particular if $x_\Omega\in X(k_\Omega)$ is a point such that $m_X(x_\Omega)=0$,
and $y_\Omega=\varphi(x_\Omega)\in Y(k_\Omega)$, then $m_Y(y_\Omega)=0$.

Let $x_\Omega\in X(k_\Omega)$ be a point,
let $S$ be a finite set of nonarchimedean places of $k$, and let
$\sU_{X,S}$ be an open neighbourhood of the $S$-part  $x_S$ of $x_\Omega$.
For $v\in\Omega_r$ we
denote by $\sU_{X,v}$ the connected component of $x_v$. We set
$\sU_{X,r}=\prod_{v\in\Omega_r} \sU_{X,v}$. We set
$\Sigma=S\cup\Omega_r$ and
$$
\sU_{X,\Sigma}=\sU_{X,S}\times \sU_{X,r}\subset X(k_\Sigma).
$$
Then $\sU_{X,\Sigma}$ is an open neighbourhood of $x_\Sigma$.
We  say that $\sU_{X,\Sigma}$ is
 \emph{the special neighbourhood of $x_\Sigma$ defined by} $\sU_{X,S}$.

For the sake of the argument it will be convenient to introduce
Property (P):

(P) \emph{For any point $x_\Omega\in X(k_\Omega)$  such that $m_X(x_\Omega)=0$,
for any finite set $S$ of nonarchimedean places of $k$, and for any open
neighbourhood $\sU_{X,S}$ of $x_S$, there exists a $k$-point
$x_0\in X(k)\cap\sU_{X,\Sigma}$, where $\sU_{X,\Sigma}$ is the special
neighbourhood of $x_\Sigma$ defined by $\sU_{X,S}$.}

  Theorem \ref{thm:Harari-generalized} precisely says that property  (P)  holds for any
 $X$ as in the theorem.

 We need a few lemmas.

\begin{lemma}\label{lem:connected-components}
Let $\psi\colon G\to G'$ be a surjective homomorphism
of $\Rr$-groups.
Let $X$ be a  homogeneous $G$-variety, and $X'$ a homogeneous $G'$-variety.
 Let $\varphi\colon X\to X'$ be a $\psi$-equivariant morphism.
Let $x\in X(\Rr)$ and set $x'=\varphi(x)\in X'(\Rr)$.
Then $\varphi$ takes the connected component of $x$ in $X(\Rr)$
onto the connected component of $x'$ in $X'(\Rr)$.
\end{lemma}
\medskip

\emph{Proof\/}\ Consider the morphism $\lambda_x\colon G\to X$
defined by $g\mapsto xg$ for $g\in G$. The morphism  $\lambda_x$
is smooth, hence the map $\lambda_x\colon G(\Rr)\to X(\Rr)$ is
open. We see that the orbit $x G(\Rr)^0$ is open,
where $G(\Rr)^0$ is the connected component of 1 in $G(\Rr)$.
Clearly $x G(\Rr)^0$ is connected.
Since all the other orbits of $G(\Rr)^0$ are also open,
we see that our orbit $x G(\Rr)^0$ is closed, hence it is the
connected component of $x$ in $X(\Rr)$. We have proved that the
connected components in $X(\Rr)$ are orbits of $G(\Rr)^0$.
Similarly the connected components in $X'(\Rr)$ are orbits of
$G'(\Rr)^0$.
Consider the action of $G$ on $G'$ by $g'\cdot g=g'\psi(g)$,
where $g'\in G',\ g\in G$.
By what has been proved, the connected component
$G'(\Rr)^0$ of 1 in $G'(\Rr)$ is the $G(\Rr)^0$-orbit of 1.
Thus $G'(\Rr)^0=\psi(G(\Rr)^0)$. Together with
the  formula $x'\psi(g)=\varphi(xg)$
this shows that $\varphi$ maps  a $G(\Rr)^0$-orbit  in $X(\Rr)$
onto a  $G'(\Rr)^0$-orbit in $X'(\Rr)$.
Thus $\varphi$ maps the connected component of  $x \in X(\Rr)$
onto the  connected component  of $\varphi(x) \in X'(\Rr)$.
QED.
\bigskip

The following lemma goes back to Cassels and Tate.

\begin{lemma}\label{lem:abelian-varieties}
Let $\psi\colon A\to A'$ be a surjective homomorphism
of abelian varieties over a number field $k$.
If $\Sha(A)$ is finite, then $\Sha(A')$ is also finite.
\end{lemma}

\emph{Proof\/} By Poincar\'e's complete reducibility theorem (cf.
\cite{Mumford},  Ch. IV, Sect. 19, Theorem 1, p. 173) there exists
an abelian variety $A''$ over $k$ such that $A$ is isogenous to
$A'\times A''$. By \cite{Milne}, Ch. I, Proof of Lemma 7.1, it
follows that $\Sha(A'\times A'')$ is finite. Thus $\Sha(A')$ is
finite. QED.

For the sake of completeness, let us give a proof of the following
well known result.

\begin{lemma}\label{lem:torsor-Brauer-one}
Let $\varphi\colon Z\to X$ be a torsor under a
quasi-trivial torus $P$,
where $Z$ and $X$ are smooth $k$-varieties
over a field $k$ of characteristic $0$.
Then there is an induced  homomorphism
$\varphi^*\colon \Bro(X^c)\to\Bro(Z^c)$, and that homomorphism
is an isomorphism.
\end{lemma}

\medskip
{\it Proof\/} Let $Y$ be a dense open set of a smooth, proper,
geometrically integral variety $Y^c$.
Let $k(Y)$ be the function field of $Y$. By well known results of
Grothendieck \cite{Groth2}
the morphisms ${\rm Spec\,} k(Y) \to Y  \to Y^{c}$ induce injections
$\Br Y^{c} \subset  \Br Y \subset \Br k(Y)$. Let $Z,Z^c,X,X^c$ be as
above. By properness of $X^c$ and smoothness of $Z ^c$,
the projection morphism $\varphi : Z \to X$ extends to a   morphism
$\varphi : W  \to X^c$, where $W \subset Z^c$ is an open set which
contains all points of codimension 1  of $Z^c$.
We thus have a natural map $\varphi^* : \Br X^c \to \Br Z$.
By the purity theorem for the Brauer group
(\cite{Groth2}, see  \cite{CTBIPG}, Section 3.4)  the restriction map $\Br Z^c
\to  \Br W$ is an isomorphism. We thus have a homomorphism
$\varphi^* : \Br  X^c \to \Br Z^c$. The map $\varphi^* : \Br  X^c \to
\Br Z^c$ is induced by the map $\Br k(X) \to \Br k(Z)$. It is none
other than the natural map of unramified cohomology groups ${\rm Br}_{\nr}
(k(X)/k)   \to {\rm Br}_{\nr}(k(Z)/k)$ (see \cite{CTBIPG}, Sections 2.2.1 and
2.2.2  and Prop. 4.2.3 (a)).

Since  $P$ is a quasi-trivial torus, over any field $F$ containing $k$,
Shapiro's lemma and Hilbert's Theorem 90 yield $\H^1_{\et}(k(X),P)=0$.
The generic fibre of $Z \to X$ is thus $k(X)$-isomorphic to
$P \times_{k}k(X)$. Since the quasi-trivial torus $P$ as a $k$-variety
is an open set of affine space over $k$, we see that
the field extension $k(Z)/k(X)$ is purely transcendental.  From
  Theorem 4.1.5  of  \cite{CTBIPG}  we get that the map $\varphi^* :
{\rm Br}_{\nr}k(X)  \to {\rm Br}_{\nr}k(Z)$ is an isomorphism. Thus   $\varphi^* : \Br  X^c \to \Br Z^c$ is an isomorphism (use \cite{CTBIPG},
Prop. 4.2.3 (a)).
Similarly $\varphi^* : \Br  {\overline X}^c \to \Br {\overline Z}^c$
is an isomorphism.
By the very definition of $\Bro$, we conclude that there is an induced map
$\varphi^* : \Bro X^c  \to \Bro Z^c $ and that this map is an isomorphism.
QED.

\medskip

We start proving Theorem \ref{thm:Harari-generalized}.
The proof is similar to that of Theorem \ref{t2}.

\medskip

\emph{First reduction.}

Let $X$ and $G$ be as in the theorem.
We write $G\uu$ for $L\uu$, where $L=G\lin$.
Set $G'=G/G\uu$, $Y=X/G\uu$.
We have a canonical smooth morphism $\varphi\colon X\to Y$.
Then $Y$ is a homogeneous space of $G'$ with connected geometrical stabilizers.
We have  $(G')\lin=G\lin/G\uu$,
hence $(G')\uu=1$.
We have $(G')\ab=G\ab$, hence $\Sha((G')\ab)$ is finite.

Assume that $Y$ has Property (P).
We prove that $X$ has this property.
Let $x_\Omega\in X(k_\Omega)$ be a point such that $m_X(x_\Omega)=0$.
Set $y_\Omega=\varphi(x_\Omega)\in Y(k_\Omega)$.
Since $m_X(x_\Omega)=0$, we see that $m_Y(y_\Omega)=0$.
Let $S$ and $\sU_{X,S}$ be as in (P).
Set $\sU_{Y,S}=\varphi(\sU_{X,S})\subset Y(k_S)$.
Since the morphism $\varphi\colon X\to Y$ is smooth,
the map $\varphi\colon X(k_S)\to Y(k_S)$ is open,
hence $\sU_{Y,S}$ is open in $Y(k_S)$.
Set $\Sigma=S\cup\Omega_r$.
Let $\sU_{Y,\Sigma}$ denote the special open neighbourhood of $y_\Sigma$
defined by $\sU_{Y,S}$.
For each $v\in\Omega_r$ let $\sU_{X,v}$ denote the connected component
of $x_v$ in $X(k_v)$.
By Lemma \ref{lem:connected-components}, for each
$v\in\Omega_r$ the set $\varphi(\sU_{X,v})$ is the
connected component of $y_v$ in $Y(k_v)$.
Thus $\sU_{Y,\Sigma}=\varphi(\sU_{X,\Sigma})$.
Since $Y$ has Property (P),
there exists a $k$-point $y_0\in Y(k)\cap\sU_{Y,\Sigma}$.

Let $X_{y_0}$ denote the fibre of $X$ over $y_0$. It is a
homogeneous space of the unipotent group  $G\uu$.
By \cite{B96}, Lemma 3.1, the $k$-variety $X_{y_0}$ has a $k$-point
and has the weak approximation property.
Consider the set $\sV_\Sigma:=
X_{y_0}(k_\Sigma)\cap\sU_{X,\Sigma}$, it is open in
$X_{y_0}(k_\Sigma)$. Since $y_0\in\varphi(\sU_{X,\Sigma})$, the
set $\sV_\Sigma$ is non-empty. Since $X_{y_0}$ has the weak
approximation property, there is a point $x_0\in X_{y_0}(k)\cap\sV_\Sigma$.
Clearly $x_0\in X(k)\cap\sU_{X,\Sigma}$.
Thus $X$ has Property (P).
Thus in the proof of Theorem \ref{thm:Harari-generalized}
we may assume    $G\uu=1$.
\medskip

\emph{Second reduction.}

By Proposition \ref{p1} we may regard $X$ as a homogeneous space
of another connected linear group $G'$ such that $(G')^\der$
is semisimple simply connected, and the stabilizers of
the geometric points of $X$ in $G'$ are linear and connected.
It follows from the construction in the proof of Proposition \ref{p1}
that there is a surjective homomorphism $G\ab\to (G')\ab$.
Since by assumption $\Sha(G\ab)$ is finite, we obtain from Lemma
\ref{lem:abelian-varieties} that $\Sha((G')\ab)$ is finite.
Thus if Theorem \ref{thm:Harari-generalized} holds for the pair $(G',X)$,
then it holds for $(G,X)$.
We see that  in the proof of Theorem \ref{thm:Harari-generalized}
we may assume that $G\lin$ is reductive, $G^\der$ is semisimple simply connected,
and the stabilizers of the geometric points of $X$ in $G$ are linear and connected.

\medskip
{\it Relaxing the assumptions}

To prove Theorem  \ref{thm:Harari-generalized} it is  enough to prove
the following result.
We write $G\sss$ for $L\sss$,  where $L=G\lin$.
The notation $\Hbar_1$ was defined in Subsection \ref{sec:structure}

\begin{theorem} \label{thm:non-connected-Manin}
Let $k$ be a number field,
$G$ a connected $k$-group, and $X$ a homogeneous space of $G$ with
geometric stabilizer $\Hbar$. Assume

{\rm (i)} $G\uu=\{1\}$,

{\rm (ii)} $\Hbar \subset \ov G\lin$,

{\rm (iii)} $G\sss$ is simply connected,

{\rm (iv)} $\Hbar_1$ is connected and
has no non-trivial characters (e.g. $\Hbar$ is connected).

{\rm (v)} $\Sha(G\ab)$ is finite.

\noindent Then $X$ has Property (P).
\end{theorem}

Recall that the homogeneous space $X$ defines a $k$-form of $\Hbar^{\rm mult}$
which we denote by $M$ (see \cite{B96}, 4.1) and that there is a natural
 homomorphism $M\to G\sab$.
We first prove a special case of Theorem \ref{thm:non-connected-Manin}.

\begin{proposition} \label{5.4-Manin} With the hypotheses of Theorem
$\ref{thm:non-connected-Manin}$, assume that $M$ {\bf injects} into
$G\sab$ (i.e. $\Hbar\cap {\ov G}\sss=\Hbar_1$).
Then $X$ has Property (P).
\end{proposition}

{\em Proof\/} Set $Y=X/G\sss$. Then $Y$ is a homogeneous space of
the semiabelian variety $G\sa$, hence it is a torsor of some
semiabelian variety $G'$. We have $(G')\ab=G\ab$, hence
$\Sha((G')\ab)$ is finite.
We have a canonical smooth morphism
$\varphi\colon X\to Y$.

Let $x_\Omega\in X(k_\Omega)$ be a point such that $m_X(x_\Omega)=0$.
Let $S$, $\sU_{X,S}$ and $\sU_{X,\Sigma}$ be as in (P).
Set $y_\Omega=\varphi(x_\Omega)\in Y(k_\Omega)$.
Since $m_X(x_\Omega)=0$, we see that $m_Y(y_\Omega)=0$.
As in the  first reduction, we define
$\sU_{Y,S}:=\varphi(\sU_{X,S})$, construct  the corresponding
special open neighbourhood $\sU_{Y,\Sigma}$ of $y_\Sigma$,
and prove that $\sU_{Y,\Sigma}=\varphi(\sU_{X,\Sigma})$.
Now since $Y$ is a torsor of a semiabelian variety
with finite Tate--Shafarevich group, by the theorem of Harari
\cite{H} the variety $Y$ has Property (P).
It follows that there exists a $k$-point $y_0\in Y(k)\cap\sU_{Y,\Sigma}$.

Let $X_{y_0}$ denote the fibre of $X$ over $y_0$.
Consider the set $\sV_\Sigma:=X_{y_0}(k_\Sigma)\cap\sU_{X,\Sigma}$,
it is open in $X_{y_0}(k_\Sigma)$.
Since $y_0\in\varphi(\sU_{X,\Sigma})$, the set $\sV_\Sigma$ is non-empty.
In particular, $X_{y_0}(k_v)\neq\emptyset$ for any $v\in\Omega_r$.
The variety $X_{y_0}$ is a homogeneous space of $G\sss$ with geometric stabilizer
$\Hbar\cap\Gbar\sss=\Hbar_1$. The group $G\sss$ is semisimple simply
connected by (iii). The group $\Hbar_1$ is connected and has no
nontrivial characters by (iv).
By \cite{B93},  Cor. 7.4, the fact that $X_{y_0}$ has points in all
real completions of $k$ is enough to ensure  that
 $X_{y_0}$  has a $k$-point. By \cite{B90}, Theorems 1.1
and 1.4 (see also \cite{CTX}), the variety $X_{y_0}$ has the weak approximation property,
and therefore  there is a point $x_0\in X_{y_0}(k)\cap\sV_\Sigma$.
Clearly $x_0\in X(k)\cap\sU_{X,\Sigma}$,
which shows that $X$ has the property (P). QED
\medskip

Let us resume the proof of Theorem \ref{thm:non-connected-Manin}.
We construct a quasi-trivial $k$-torus $P$,
the $k$-group $F:=G\times P$,
a homogeneous space $Z$ of $F$, and a morphism $\pi\colon Z\to X$
as in the proof of Theorem \ref{thm:non-connected}.
Since $(Z,\pi)$ is a torsor under the quasi-trivial torus $P$,
by Lemma \ref{lem:torsor-Brauer-one}
the canonical map
\begin{equation*}
\pi_*\colon \Bro (Z^c)^D\to\Bro (X^c)^D.
\end{equation*}
is an isomorphism.
We have $F\ab=G\ab$, hence $\Sha(F\ab)$ is finite.

Let $x_\Omega\in X(k_\Omega)$ be a point,
and assume that $m_X(x_\Omega)=0$.
Since $\pi\colon Z\to X$ is a torsor under a quasi-trivial torus,
we can lift $x_\Omega$ to some $z_\Omega\in Z(k_\Omega)$.
We have $m_X(x_\Omega)=\pi_*(m_Z(z_\Omega))$.
Since $\pi_*$ is an isomorphism,
from $m_X(x_\Omega)=0$ we conclude that $m_Z(z_\Omega)=0$.

Let $S$ be as above, and let $\sU_{X,S}\subset X(k_S)$
be an open neighbourhood of $x_S$.
Let $\sU_{X,\Sigma}\subset X(k_\Sigma)$
be the corresponding special neighbourhood of $x_\Sigma$.
Set
$$
\sU_{Z,S}=\pi^{-1}(\sU_{X,S})\subset Z(k_S).
$$
For $v\in\Omega_r$ let $\sU_{Z,v}$ be
the connected component of $z_v$ in $Z(k_v)$.
By Lemma \ref{lem:connected-components} $\pi(\sU_{Z,v})=\sU_{X,v}$.
Set $\sU_{Z,r}=\prod_{v\in\Omega_r} \sU_{Z,v}$ and
$\sU_{Z,\Sigma}=\sU_{Z,S}\times\sU_{Z,r}$.
Then $\sU_{Z,\Sigma}$ is a special open neighbourhood of $z_\Sigma$,
and $\pi(\sU_{Z,\Sigma})=\sU_{X,\Sigma}$.

The homogeneous space $Z$ of $F$ satisfies the hypotheses
of Proposition \ref{5.4-Manin},
so by this proposition there is a point $z_0\in Z(k)\cap\sU_{Z,\Sigma}$.
Set $x_0=\pi(z_0)$, then $x_0\in X(k)\cap\sU_{X,\Sigma}$.
Thus $X$ has Property (P).
This completes the proofs of Theorem \ref{thm:non-connected-Manin}
and Theorem \ref{thm:Harari-generalized}.
QED.

\bigskip

\noindent Raymond and Beverly Sackler School of Mathematical
Sciences, Tel Aviv University, 69978 Tel Aviv, Israel

\noindent  borovoi@post.tau.ac.il
\bigskip

\noindent
CNRS, UMR 8628, Math\'ematiques, B\^atiment 425, Universit\'e Paris-Sud,
F-91405 Orsay, France

\noindent Jean-Louis.Colliot-Thelene@math.u-psud.fr

\bigskip

\noindent Department of Mathematics, South Kensington Campus,
Imperial College London,
SW7 2BZ England, U.K.

\smallskip

\noindent Institute for the Information Transmission Problems,
Russian Academy of Sciences, 19 Bolshoi Karetnyi,
Moscow, 127994 Russia
\medskip

\noindent a.skorobogatov@imperial.ac.uk

\end{document}